\documentclass[12pt]{article}
\usepackage{amsfonts}
\usepackage{amsmath}
\usepackage{amssymb}
\usepackage{graphicx}

\providecommand{\U}[1]{\protect\rule{.1in}{.1in}}
\newtheorem{theorem}{Theorem}

\newtheorem{corollary}[theorem]{Corollary}

\newtheorem{definition}[theorem]{Definition}

\newtheorem{proposition}[theorem]{Proposition}
\newtheorem{remark}[theorem]{Remark}

\newenvironment{proof}[1][Proof]{\noindent\textbf{#1.} }{\ \rule{0.5em}{0.5em}}

\begin{document}

\title{Stochastic approach for a multivalued Dirichlet-Neumann problem}
\author{Lucian Maticiuc$^{{\scriptsize a,\ast }}$, Aurel R\u{a}\c{s}canu$^{%
{\scriptsize b,\ast }}$ \\
\bigskip \\
$^{{\scriptsize a}}${\scriptsize Department of Mathematics,
\textquotedblleft Gh. Asachi\textquotedblright\ Technical University of Ia%
\c{s}i, Bd. Carol I, no. 11, 700506, Romania,}\\
$^{{\scriptsize b}}${\scriptsize Faculty of Mathematics, \textquotedblleft
Al.I. Cuza\textquotedblright\ University of Ia\c{s}i, Bd. Carol I, no. 9,
700506, and}\\
{\scriptsize \textquotedblleft O. Mayer\textquotedblright Mathematics
Institute of the Romanian Academy, Ia\c{s}i, Bd. Carol I, no. 8, 700506,
Romania.}}
\maketitle

\footnotetext[1]{{\scriptsize The work was supported by Grant ID
395/2007, CEEX 06-11-56/2006.}}

\footnotetext{\textit{\scriptsize E-mail addresses:} {\scriptsize %
lucianmaticiuc@yahoo.com (Lucian Maticiuc),~aurel.rascanu@uaic.ro (Aurel R%
\u{a}\c{s}canu),}}

\begin{abstract}
We prove the existence and uniqueness of a viscosity solution of the
parabolic variational inequality\ (PVI)\ with a mixed nonlinear multivalued
Neumann-Dirichlet boundary condition:%
\begin{equation*}
\left\{
\begin{array}{r}
\dfrac{\partial u(t,x)}{\partial t}-\mathcal{L}_{t}u\left( t,x\right) {+}{%
\partial \varphi }\big(u(t,x)\big)\ni f\big(t,x,u(t,x),\!(\nabla u\sigma
)(t,x)\big),\;t>0,\;x\in \mathcal{D},\medskip \\
\multicolumn{1}{l}{\dfrac{\partial u(t,x)}{\partial n}+{\partial \psi }%
\big(u(t,x)\big)\ni g\big(t,x,u(t,x)\big),\;\;t>0,\;x\in Bd\left( \mathcal{D}%
\right) ,\vspace*{2mm}} \\
\multicolumn{1}{l}{u(0,x)=h(x),\;x\in \overline{\mathcal{D}},}%
\end{array}%
\right.
\end{equation*}%
where $\partial \varphi $ and $\partial \psi $ are subdifferentials
operators and $\mathcal{L}_{t}$ is a second differential operator given by%
\begin{equation*}
\mathcal{L}_{t}v\left( x\right) =\frac{1}{2}\sum_{i,j=1}^{d}(\sigma \sigma
^{\ast })_{ij}(t,x)\frac{{\partial }^{2}v\left( x\right) }{{\partial }x_{i}{%
\partial }x_{j}}+\sum_{i=1}^{d}b_{i}(t,x)\frac{{\partial }v\left( x\right) }{%
{\partial }x_{i}}.
\end{equation*}%
The result is obtained by a stochastic approach. First we study the
following backward stochastic generalized variational inequality:%
\begin{equation*}
\left\{
\begin{array}{l}
dY_{t}{+}F\left( t,Y_{t},Z_{t}\right) dt{+}G\left( t,Y_{t}\right) dA_{t}\in
\partial \varphi \left( Y_{t}\right) dt{+}\partial \psi \left( Y_{t}\right)
dA_{t}{+}Z_{t}dW_{t}~,\;0\leq t\leq T,\medskip \\
Y_{T}=\xi ,%
\end{array}%
\right.
\end{equation*}%
where $\left( A_{t}\right) _{t\geq 0}$ is a continuous one-dimensional
increasing measurable process, and then we obtain a Feynman-Ka\c{c}
representation formula for the viscosity solution of the PVI problem.
\end{abstract}

\textbf{AMS Classification subjects: }35D05, 35K85, 60H10, 60H30, 47J20,
49J40.\medskip


\section{Introduction}

Viscosity solutions were introduced by M.G. Crandall and P.L. Lions in \cite%
{cr-li/83}, and then developed in the classical work of M.G. Crandall, H. Ishii, P.L. Lions \cite{cr-li/92}, where are presented several equivalent ways to formulate the
notion of such type solutions. The framework of this theory allows for
merely continuous functions to be the solutions of fully nonlinear equations
of second order which provides a very general existence and uniqueness
theorems.

In 1992 E. Pardoux and S. Peng \cite{pa-pe/92} introduced backward
stochastic differential equations (BSDE) and supplied probabilistic formulas
for the viscosity solutions of semilinear partial differential equations,
both of parabolic and elliptic type in whole space. Elliptic equations with
Dirichlet boundary condition have been treated by R.W.R. Darling and E.
Pardoux in \cite{da-pa/97} and with a homogeneous Neumann boundary condition
by Y. Hu in \cite{hu/93}.

The parabolic (and elliptic) systems of partial differential equations (PVI
without the subdifferential operator) with nonlinear Neumann boundary
conditions was the subject of the paper E. Pardoux and S. Zhang \cite%
{pa-zh/98}. The case of the systems of variational inequalities for partial
differential equations in whole space was studied by L. Maticiuc, E.
Pardoux, A. R\u{a}\c{s}canu and A. Z\u{a}linescu in \cite{ma-pa-ra-za/08}.

The main idea for proving the existence of the viscosity solutions
for PDE and PVI is the stochastic approach. Using a suitable BSDE,
or backward stochastic variational inequality (BSVI) for the PVI
case, one can obtain a generalizations of the Feynman-Ka\c{c}
formula (i.e. stochastic representation formula of the viscosity
solution for deterministic problems).

The origin of our study goes from the PDE%
\[
\left\{
\begin{array}
[c]{r}%
\dfrac{\partial u}{\partial t}-\mathcal{L}_{t}u=f,\;t>0,\;x\in\mathcal{D}%
,\medskip\\
\multicolumn{1}{l}{\dfrac{\partial u}{\partial n}=g,\;t>0,\;x\in
Bd\left(
\mathcal{D}\right)  ,\vspace*{2mm}}\\
\multicolumn{1}{l}{u(0,x)=h(x),\;x\in\overline{\mathcal{D}},}%
\end{array}
\right.
\]
which is a mathematical model for the evolution of a state
$u(t,x)\in\mathbb{R}$ of a diffusion dynamical system with sources
$f$ acting in the interior of the domain $\mathcal{D}$ and $g$ on
the boundary of $\mathcal{D}$.

In certain applications it is call upon to maintain the state
$u(t,x)$ in a interval $\mathbb{I}\subset\mathbb{R}$ for all
$x\in\mathcal{D}$ and in a
interval $\mathbb{J}\subset\mathbb{R}$ for all $x\in Bd\left(  \mathcal{D}%
\right)  $. Practically these can be realized adding the
supplementary sources $\partial I_{\mathbb{I}}\left(  u\left(
t,x\right)  \right) $ and $\partial I_{\mathbb{J}}\left(  u\left(
t,x\right)  \right) $ on the system. These sources produce
\textquotedblleft inward pushes\textquotedblright\ that would
keep the state process%
\[
u(t,x)\text{ in }\mathbb{I}\text{, }\forall x\in\mathcal{D}\quad
\text{and}\quad u(t,x)\text{ in }\mathbb{J}\text{, }\forall x\in
Bd\left( \mathcal{D}\right)
\]
and do this in a minimal way (i.e. only when $u(t,x)$ arrives on the
boundary
of $\mathbb{I}$ and respectively $\mathbb{J}$). Hence $\partial I_{\mathbb{I}%
}\left(  u\left(  t,x\right)  \right)  $ and $\partial
I_{\mathbb{J}}\left( u\left(  t,x\right)  \right)  $ represent
perfect feedback flux controls.

The aim of this paper is to treat the more general case of a
parabolic variational inequality\ with mixed nonlinear multivalued
Neumann-Dirichlet boundary condition. This requires the presence of
a new terms in the associated BSVI considered, namely an integral
with respect to a continuous increasing process.

The scalar BSDE with one-sided reflection, which provides a probabilistic
representation for the unique viscosity solution of an obstacle problem for
a nonlinear parabolic PDE, was considered by N. El
Karoui, C. Kapoudjian, E. Pardoux, S. Peng, M.C. Quenez in \cite{ka-ka-pa/97}%
. E. Pardoux and A. R\u{a}\c{s}canu in \cite{pa-ra/98} (and \cite{pa-ra/99} for the generalization
to the Hilbert spaces framework) studied the general case of
BSVI and obtained probabilistic representation for the solutions of PVI in whole space.

The paper is organized as follows: In Section 2 we formulate the
Neumann-Dirichlet PVI problem; we present the main results and we prove the uniqueness theorem. For the existence theorem we first study in Section 3 a
certain BSVI. The solution of this backward equation gives us, via Feynman-Ka\c{c}
representation formula, a viscosity solution for the deterministic
multivalued partial differential equation as it is shown in Section 4.

\section{Main results}

\label{sect1}Let $\mathcal{D}$ be a open connected bounded subset of $%
\mathbb{R}^{d}$ of form%
\begin{equation*}
\mathcal{D}=\left\{ x\in \mathbb{R}^{d}:\ell \left( x\right) <0\right\}
,\;\;Bd\left( \mathcal{D}\right) =\left\{ x\in \mathbb{R}^{d}:\ell \left(
x\right) =0\right\} ,
\end{equation*}%
where $\ell \in C_{b}^{3}\left( \mathbb{R}^{d}\right) $, $\left\vert \nabla
\ell \left( x\right) \right\vert =1,\;$for all $x\in Bd\left( \mathcal{D}%
\right) $.

We define outward normal derivative by%
\begin{equation*}
\frac{\partial v\left( x\right) }{\partial n}=\sum_{j=1}^{d}\frac{\partial
\ell \left( x\right) }{\partial x_{j}}\frac{\partial v\left( x\right) }{%
\partial x_{j}}=\left\langle \nabla \ell \left( x\right) ,\nabla v\left(
x\right) \right\rangle ,\;\text{for all }x\in Bd\left( \mathcal{D}\right) .
\end{equation*}%
The aim of this paper is to study the existence and uniqueness of a
viscosity solution for the following parabolic variational inequality (PVI)
with a mixed nonlinear multivalued Neumann-Dirichlet boundary condition:$%
\smallskip $%
\begin{equation}
\left\{
\begin{array}{r}
\dfrac{\partial u(t,x)}{\partial t}-\mathcal{L}_{t}u\left( t,x\right) +{%
\partial \varphi }\big(u(t,x)\big)\ni f\big(t,x,u(t,x),(\nabla u\sigma )(t,x)%
\big), \\
t>0,\;x\in \mathcal{D},\medskip \\
\multicolumn{1}{l}{\dfrac{\partial u(t,x)}{\partial n}+{\partial \psi }\big(%
u(t,x)\big)\ni g\big(t,x,u(t,x)\big),\;\;t>0,\;x\in Bd\left( \mathcal{D}%
\right) ,\vspace*{2mm}} \\
\multicolumn{1}{l}{u(0,x)=h(x),\;x\in \overline{\mathcal{D}},}%
\end{array}%
\right.  \label{pvi}
\end{equation}%
where operator $\mathcal{L}_{t}$ is given by%
\begin{align*}
\mathcal{L}_{t}v(x)& =\dfrac{1}{2}\mathrm{Tr}\big[\sigma (t,x)\sigma ^{\ast
}(t,x)D^{2}v(x)\big]+\big\langle b(t,x),\nabla v(x)\big\rangle\medskip \\
& =\dfrac{1}{2}\sum\limits_{i,j=1}^{d}(\sigma \sigma ^{\ast })_{ij}(t,x)%
\dfrac{\partial ^{2}v(x)}{\partial x_{i}\partial x_{j}}+\sum%
\limits_{i=1}^{d}b_{i}(t,x)\dfrac{\partial v(x)}{\partial x_{i}}\,.
\end{align*}%
for $v\in C^{2}\left( \mathbb{R}^{d}\right) .$

We will make the following assumptions:

\begin{itemize}
\item[(I)] \ Functions%
\begin{equation}
\begin{array}{l}
b:\left[ 0,\infty \right) \times \mathbb{R}^{d}\rightarrow \mathbb{R}%
^{d},\medskip \\
\sigma :\left[ 0,\infty \right) \times \mathbb{R}^{d}\rightarrow \mathbb{R}%
^{d\times d},\medskip \\
f:\left[ 0,\infty \right) \times \overline{\mathcal{D}}\times \mathbb{R}%
\times \mathbb{R}^{d}\rightarrow \mathbb{R},\medskip \\
g:\left[ 0,\infty \right) \times Bd\left( \mathcal{D}\right) \times \mathbb{R%
}\rightarrow \mathbb{R},\medskip \\
h:\overline{\mathcal{D}}\rightarrow \mathbb{R\;\;\;\;}\text{are continuous}%
\end{array}
\label{h1}
\end{equation}%
We assume that for all $T>0$ there exist $\alpha \in \mathbb{R}$ and $%
L,\beta ,\gamma \geq 0$ (which can depend on $T$) such that $\forall t\in %
\left[ 0,T\right] ,\;\forall x,\tilde{x}\in \mathbb{R}^{d}:$%
\begin{equation}
\big|b\left( t,x\right) -b\left( t,\tilde{x}\right) \big|+\big\|\sigma
\left( t,x\right) -\sigma \left( t,\tilde{x}\right) \big\|\leq L\left\Vert x-%
\tilde{x}\right\Vert ,  \label{h2}
\end{equation}%
and$\ \forall t\in \left[ 0,T\right] $, $\forall x\in \overline{\mathcal{D}}$%
, $u\in Bd\left( \mathcal{D}\right) $, $y,\tilde{y}\in \mathbb{R},z,\tilde{z}%
\in \mathbb{R}^{d}$:%
\begin{equation}
\begin{array}{rl}
\left( i\right) \text{\ \ } & (y-\tilde{y})\big(f(t,x,y,z)-f(t,x,\tilde{y},z)%
\big)\leq \alpha |y-\tilde{y}|^{2},\vspace*{2mm} \\
\left( ii\right) \text{\ \ } & \big|f(t,x,y,z)-f(t,x,y,\tilde{z})\big|\leq
\beta |z-\tilde{z}|,\vspace*{2mm} \\
\left( iii\right) \text{\ \ } & \big|f(t,x,y,0)\big|\leq \gamma \big (1+|y|%
\big ),\vspace*{2mm} \\
\left( iv\right) \text{\ \ } & (y-\tilde{y})\big(g(t,u,y)-g(t,u,\tilde{y})%
\big)\leq \alpha |y-\tilde{y}|^{2},\vspace*{2mm} \\
\left( v\right) \text{\ \ } & \big|g(t,u,y)\big|\leq \gamma \big(1+|y|\big).%
\vspace*{2mm}%
\end{array}
\label{h3}
\end{equation}%
In fact, condition (\ref{h3}-$i$) and (\ref{h3}-$iv$) mean that, for all $%
t\in \left[ 0,T\right] $, $x\in \overline{\mathcal{D}}$, $u\in Bd\left(
\mathcal{D}\right) $, $z\in \mathbb{R}^{d}$:%
\begin{equation*}
\begin{array}{l}
y\mapsto \alpha y-f\left( t,x,y,z\right) :\mathbb{R\rightarrow R}\medskip \\
y\mapsto \alpha y-g\left( t,u,y\right) :\mathbb{R\rightarrow R}%
\end{array}%
\end{equation*}%
are increasing functions.

\item[(II)] With respect to functions $\varphi $ and $\psi $ we assume%
\begin{equation}
\begin{array}{rl}
\left( i\right) \text{\ \ } & \varphi ,\psi :\mathbb{R}\rightarrow (-\infty
,+\infty ]\text{ \ are proper convex l.s.c. functions,}\medskip \\
\left( ii\right) \text{\ \ } & \varphi \left( y\right) \geq \varphi \left(
0\right) =0\text{ and }\psi \left( y\right) \geq \psi \left( 0\right) =0,\
\forall \;y\in \mathbb{R},\medskip%
\end{array}
\label{h4}
\end{equation}%
and there exists a positive constant $M$ such that%
\begin{equation}
\begin{array}{rl}
\left( i\right) \text{\ \ } & \Big|\varphi \big (h(x)\big )\Big|\leq
M,\;\;\forall {x}\in \overline{\mathcal{D}},\medskip \\
\left( ii\right) \text{\ \ } & \Big|\psi \big (h(x)\big )\Big|\leq
M,\;\;\forall {x}\in Bd\left( \mathcal{D}\right) .\medskip%
\end{array}
\label{h5}
\end{equation}
\end{itemize}

\begin{remark}
\label{obs1}Condition (\ref{h4}-ii) is generally realized by changing
problem (\ref{pvi}) in an equivalent form, as example: if $\left(
u_{0},u_{0}^{\ast }\right) \in \partial \varphi $ we can replace $\varphi
\left( u\right) $ by $\varphi \left( u+u_{0}\right) -\varphi \left(
u_{0}\right) -\left\langle u_{0}^{\ast },u\right\rangle $; a similar
transformation one can do for $\psi $.
\end{remark}

We denote%
\begin{equation*}
\begin{array}{l}
Dom\left( \varphi \right) =\left\{ u\in \mathbb{R}:\varphi \left( u\right)
<\infty \right\} ,\medskip \\
\partial \varphi \left( u\right) =\left\{ u^{\ast }\in \mathbb{R}:u^{\ast
}\left( v-u\right) +\varphi \left( u\right) \leq \varphi \left( v\right)
,\forall v\in \mathbb{R}\right\} ,\medskip \\
Dom\left( \partial \varphi \right) =\left\{ u\in \mathbb{R}:\partial \varphi
\left( u\right) \neq \emptyset \right\} ,\medskip \\
\left( u,u^{\ast }\right) \in \partial \varphi \Leftrightarrow u\in \
Dom\partial \varphi ,\;\;u^{\ast }\in \partial \varphi \left( u\right)
\medskip%
\end{array}%
\end{equation*}%
(for function $\psi $ we have the similar notations).

In every point $y\in Dom\left( \varphi \right) $ we have%
\begin{equation*}
{\partial \varphi }(y)=\mathbb{R}\cap \big[\varphi _{-}^{\prime }(y),\varphi
_{+}^{\prime }(y)\big],
\end{equation*}%
where $\varphi _{-}^{\prime }(y)$ and $\varphi _{+}^{\prime }(y)$ are left
derivative and, respectively, right derivative at point $y$.

\begin{itemize}
\item[(III)] We introduce the \textit{compatibility assumptions }: \newline
for all \ $\varepsilon >0$, $t\geq 0$, $x\in Bd\left( \mathcal{D}\right) $, $%
\tilde{x}\in \overline{\mathcal{D}}$, $y\in \mathbb{R}$ and $z\in \mathbb{R}%
^{d}$%
\begin{equation}
\begin{array}{cl}
\left( i\right) & \nabla \varphi _{\varepsilon }\left( y\right) g\left(
t,x,y\right) \leq \big[\nabla \psi _{\varepsilon }\left( y\right) g\left(
t,x,y\right) \big]^{+},\medskip \\
\left( ii\right) & \nabla \psi _{\varepsilon }\left( y\right) f\left( t,%
\tilde{x},y,z\right) \leq \big[\nabla \varphi _{\varepsilon }\left( y\right)
f\left( t,\tilde{x},y,z\right) \big]^{+},%
\end{array}
\label{h6}
\end{equation}%
where $a^{+}=\max \left\{ 0,a\right\} $ and $\nabla \varphi _{\varepsilon
}\left( y\right) $, $\nabla \psi _{\varepsilon }\left( y\right) $ are unique
solutions $U$ and $V$, respectively, of equations%
\begin{equation*}
{\partial \varphi }(y-\varepsilon U)\ni U\;\;\;\text{and \ \ \ }{\partial
\psi }(y-\varepsilon V)\ni V.
\end{equation*}
\end{itemize}

\begin{remark}
A) Clearly, using the monotonicity of $\nabla \varphi _{\varepsilon },\nabla
\psi _{\varepsilon }$, we see that, if
\begin{equation*}
y\cdot g\left( t,x,y\right) \leq 0\;\;\text{and}\ \ y\cdot f\left( t,\tilde{x%
},y,z\right) \leq 0
\end{equation*}%
for all $t\geq 0$, $x\in Bd\left( \mathcal{D}\right) $, $\tilde{x}\in
\overline{\mathcal{D}}$, $y\in \mathbb{R}$ and $z\in \mathbb{R}^{d}$, then
compatibility assumptions (\ref{h6}) are satisfied.$\medskip $

B) If $\varphi ,\psi :\mathbb{R}\rightarrow (-\infty ,+\infty ]$ are convex
indicator functions%
\begin{equation*}
\varphi \left( y\right) =I_{[a,\infty )}\left( y\right) =\left\{
\begin{array}{l}
0\;\;~~,\text{\ if }y\in \lbrack a,\infty )\medskip \\
+\infty ,\text{\ if }y\notin \lbrack a,\infty )%
\end{array}%
\right.
\end{equation*}%
and%
\begin{equation*}
\psi \left( y\right) =I_{(-\infty ,b]}\left( y\right) =\left\{
\begin{array}{l}
0\;\;~~,\text{\ if }y\in (-\infty ,b]\medskip \\
+\infty ,\text{\ if }y\notin (-\infty ,b]%
\end{array}%
\right.
\end{equation*}%
where $a\leq 0\leq b,$\newline
\textit{then }%
\begin{equation*}
\nabla \varphi _{\varepsilon }\left( y\right) =-\frac{1}{\varepsilon }\left(
y-a\right) ^{-}\mathit{\quad }\text{\textit{and}}\mathit{\quad }\nabla \psi
_{\varepsilon }\left( y\right) =\frac{1}{\varepsilon }\left( y-b\right) ^{+}
\end{equation*}%
and \textit{the compatibility assumptions become}%
\begin{align*}
g\left( t,x,y\right) & \geq 0,\mathit{\;}\quad \text{\textit{for} }y\leq a,\;%
\text{\textit{and\ }}\medskip \\
f\left( t,\tilde{x},y,z\right) & \leq 0,\mathit{\;}\quad \text{\textit{for} }%
y\geq b.
\end{align*}
\end{remark}

We shall define now the notion of viscosity solution in the language of sub-
and super-jets, see \cite{cr-li/92}. S\negthinspace$\mathbb{R}^{d\times d}$
will denote below the set of $d\times d$ symmetric non--negative real
matrices.

\begin{definition}
\label{jet}Let $u:[0,\infty )\times \overline{\mathcal{D}}\rightarrow
\mathbb{R}$ a continuous function, and $(t,x)\in \lbrack 0,\infty )\times
\overline{\mathcal{D}}.$ We denote by $\mathcal{P}^{2,+}u(t,x)$ \ (the
parabolic superjet of $u$ at $(t,x)$) the set of triples $(p,q,X)\in \mathbb{%
R}\times \mathbb{R}^{d}\times $S$\!\mathbb{R}^{d\times d}$ which are such
that for all $\left( s,y\right) \in \lbrack 0,\infty )\times \overline{%
\mathcal{D}}$ in a neighbourhood of $\left( t,x\right) $:%
\begin{equation*}
\begin{array}{l}
u(s,y)\leq u(t,x)+p(s-t)+\left\langle q,y-x\right\rangle +\vspace*{2mm} \\
\text{ \ \ \ \ \ \ \ \ \ \ \ \ \ \ \ }+\dfrac{1}{2}\big\langle X(y-x),y-x%
\big\rangle+o\big(|s-t|+|y-x|^{2}\big).%
\end{array}%
\end{equation*}%
Similarly is defined $\mathcal{P}^{2,-}u(t,x)$ \ (the parabolic subjet of $u$
at $(t,x)$) as the set of triples $(p,q,X)\in \mathbb{R}\times \mathbb{R}%
^{d}\times $S$\!\mathbb{R}^{d\times d}$ which are such that for all $\left(
s,y\right) \in \lbrack 0,\infty )\times \overline{\mathcal{D}}$ in a
neighbourhood of $\left( t,x\right) $:%
\begin{equation*}
\begin{array}{l}
u(s,y)\geq u(t,x)+p(s-t)+\left\langle q,y-x\right\rangle +\vspace*{2mm} \\
\;\ \ \ \ \ \ \ \ \ \ \ \ \ \ +\dfrac{1}{2}\big\langle X(y-x),y-x\big\rangle%
+o\big(|s-t|+|y-x|^{2}\big),\vspace*{2mm}%
\end{array}%
\end{equation*}%
where $r\rightarrow o\left( r\right) $ is the Landau function i.e. $o:\left[
0,\infty \right[ \rightarrow \mathbb{R}$ is a continuous function such that $%
\underset{r\rightarrow 0}{\lim }\dfrac{o\left( r\right) }{r}=0$
\end{definition}

We can give now the definition of a viscosity solution of the parabolic
variational inequality (\ref{pvi}). We denote first%
\begin{align*}
& V\left( t,x,p,q,X\right) \overset{def}{=}p-\frac{1}{2}\mathrm{Tr}\big(%
(\sigma\sigma^{\ast})(t,x)X\big)-\big\langle b(t,x),q\big\rangle \\
& \;\;\;\;\;\;\;\;\;\;\;\;\;\;\;\;\;\;\;\;\;\;\;\;\;\;\;\;\;-f\big(%
t,x,u(t,x),q\sigma (t,x)\big).
\end{align*}

\begin{definition}
Let $u:\left[ 0,\infty \right) \times \overline{\mathcal{D}}\rightarrow
\mathbb{R}$ a continuous function, which satisfies $u(0,x)=h\left( x\right)
,\;\forall ~x\in \overline{\mathcal{D}}$,\newline
(a) $u$ is a viscosity subsolution of (\ref{pvi}) if:%
\begin{equation*}
\left\vert \text{\ }%
\begin{array}{l}
u(t,x)\in Dom\left( \varphi \right) ,\ \ \forall {(t,x)}\in (0,\infty
)\times \overline{\mathcal{D}}, \\
u(t,x)\in Dom\left( \psi \right) ,\ \ \ \forall {(t,x)}\in (0,\infty )\times
Bd\left( \mathcal{D}\right) ,%
\end{array}%
\right. \vspace*{2mm}
\end{equation*}%
and, at any point $\left( t,x\right) \in (0,\infty )\times \overline{%
\mathcal{D}}$, for any $(p,q,X)\in \mathcal{P}^{2,+}u(t,x)$:%
\begin{equation}
\left\vert \text{\ }%
\begin{array}{l}
V\left( t,x,p,q,X\right) +\varphi _{-}^{\prime }\big(u(t,x)\big)\leq 0\
\;if\;x\in \mathcal{D},\medskip \\
\min \Big\{V\left( t,x,p,q,X\right) +\varphi _{-}^{\prime }\big(u(t,x)\big)%
,\;\big\langle\nabla \ell \left( x\right) ,q\big\rangle-g\big(t,x,u(t,x)\big)%
\medskip \\
\;\ \ \ \ \ \ \ \ \ \ \ \ \ \ \ \ \ \ \ \ \ \ \ \ \ \ \ \ \ \ \ \ \ \ \ \ \
\ \ +\psi _{-}^{\prime }\big(u(t,x)\big)\Big\}\leq 0\;\ if\;x\in Bd\left(
\mathcal{D}\right) .%
\end{array}%
\right.  \label{subs}
\end{equation}%
\textit{(b) \ the viscosity supersolution of (\ref{pvi}) is defined in a
similar manner as above, with }$\mathcal{P}^{2,+}$ replaced by $\mathcal{P}%
^{2,-}$, \textit{the left derivative replaced by the right derivative}, $\
\min $ \ by $\ \max $, \textit{and the inequalities} \textit{\ }$\leq \quad $%
\textit{by}$\quad \geq \ $.\newline
(c) \ a continuous function $u:\left[ 0,\infty \right) \times \overline{%
\mathcal{D}}$ is a viscosity solution of (\ref{pvi}) if it is both a
viscosity sub- and super-solution.\bigskip
\end{definition}

We now present the main results\medskip

\begin{theorem}[Existence]
\label{t2}\textit{Let assumptions (\ref{h1})-(\ref{h6}) be satisfied. Then
PVI (\ref{pvi}) has a viscosity solution.}
\end{theorem}

For the proof of the existence we shall study a certain backward stochastic
generalized variational inequality (then we use a nonlinear representation
Feynman-Ka\c{c} type formula). We present this approach in the following
section and after then the proof of Theorem \ref{t2} in Section \ref{sect3}.

\begin{theorem}[Uniqueness]
\label{t3}Let the assumptions of Theorem \ref{t2} be satisfied. If \textit{%
function}%
\begin{equation}
r\rightarrow g(t,x,r)\text{ is decreasing for }t\geq 0\text{, }x\in Bd\left(
\mathcal{D}\right) ,  \label{h8}
\end{equation}%
and \textit{there exists a continuous function }$\mathbf{m}:[0,\infty
)\rightarrow \lbrack 0,\infty )$,\textit{\ }$\mathbf{m}\left( 0\right) =0$,%
\textit{\ such that}%
\begin{equation}
\begin{array}{l}
\big|f(t,x,r,p)-f(t,y,r,p)\big|\leq \mathbf{m}\big(\left\vert x-y\right\vert
\left( 1+\left\vert p\right\vert \right) \big),\vspace*{2mm} \\
\text{\ \ \ \ \ \ }\;\text{\ \ \ \ \ \ \ }\;\text{\ \ \ \ \ \ \ }\;\text{\ \
\ \ \ \ \ }\;\text{\ }\forall \ t\geq 0,\;x,y\in \overline{\mathcal{D}}%
,\;p\in \mathbb{R}^{d},%
\end{array}
\label{h7}
\end{equation}%
\textit{then the viscosity solution is unique.}
\end{theorem}

\begin{proof}
It is sufficient to prove the uniqueness on a fixed arbitrary interval $%
\left[ 0,T\right] .$

Also, it suffices to prove that if $u$ is a subsolution and $v$ is a
supersolution such that $u(0,x)=v(0,x)=h\left( x\right) $, $x\in \overline{%
\mathcal{D}}$, then $u\leq v$.

Firstly, from definition of $\mathcal{D}$, there exists a function $\tilde{%
\ell}\in C_{b}^{3}\left( \mathbb{R}^{d}\right) $ such that $\tilde{\ell}%
\left( x\right) \geq0$ on $\overline{\mathcal{D}}$ with $\nabla\tilde{\ell}%
\left( x\right) =\nabla\ell\left( x\right) $ for $x\in Bd\left( \mathcal{D}%
\right) $ (as example $\tilde{\ell}\left( x\right) =\ell\left( x\right) +%
\underset{y\in\overline{\mathcal{D}}}{\sup}\left\vert \ell\left( y\right)
\right\vert $).

For $\lambda =\left\vert \alpha \right\vert +1$ and $\delta ,\varepsilon
,c>0 $ let
\begin{equation*}
\begin{array}{l}
\bar{u}\left( t,x\right) =e^{\lambda t}u\left( t,x\right) -\delta \tilde{\ell%
}\left( x\right) -c\vspace*{2mm} \\
\bar{v}\left( t,x\right) =e^{\lambda t}v\left( t,x\right) +\delta \tilde{\ell%
}\left( x\right) +c+\varepsilon /t.%
\end{array}%
\end{equation*}%
Denote%
\begin{equation}
\begin{array}{l}
\tilde{f}\left( t,x,r,q,X\right) =\lambda r-\dfrac{1}{2}\mathrm{Tr}\big[%
\left( \sigma \sigma ^{\ast }\right) (t,x)X\big]-\big\langle b(t,x),q%
\big\rangle\vspace*{2mm} \\
\text{\ \ \ \ \ \ }\;\text{\ \ \ \ \ \ \ }\;\text{\ \ \ \ \ \ }-e^{\lambda
t}f\big (t,x,e^{-\lambda t}r,e^{-\lambda t}q\sigma \left( t,x\right) \big )%
\end{array}
\label{defft}
\end{equation}%
and%
\begin{equation*}
\tilde{g}(t,x,r)=e^{\lambda t}g(t,x,e^{-\lambda t}r)
\end{equation*}%
Clearly $r\rightarrow \tilde{f}\left( t,x,r,q,X\right) $ is an increasing
function for all $\left( t,x,q,X\right) \in \left[ 0,T\right] \times \mathbb{%
R}^{d}\times \mathbb{R}^{d}\times $S$\!\mathbb{R}^{d\times d}$. Moreover,
since%
\begin{equation*}
M=\sup_{\left( t,x\right) \in \left[ 0,T\right] \times \overline{\mathcal{D}}%
}\left\{ |\tilde{\ell}\left( x\right) |+|D\tilde{\ell}\left( x\right)
|+|D^{2}\tilde{\ell}\left( x\right) |+\left\vert b(t,x)\right\vert
+\left\vert \sigma (t,x)\right\vert \right\} <\infty ,
\end{equation*}%
then we can choose $c=c\left( \delta ,M\right) >0$ such that for $\bar{u}=%
\bar{u}\left( t,x\right) $ and $\tilde{\ell}=\tilde{\ell}\left( x\right) :$
\begin{equation*}
\tilde{f}\left( t,x,\bar{u},D\bar{u},D^{2}\bar{u}\right) \leq \tilde{f}(t,x,%
\bar{u}+\delta \tilde{\ell}+c,D\bar{u}+\delta D\tilde{\ell},D^{2}\bar{u}%
+\delta D^{2}\tilde{\ell}).
\end{equation*}%
Using these properties$,$ assumption (\ref{h8}), and the fact that left and
right derivative of $\varphi $,$\psi $ are increasing we infer that function
$\bar{u}$ satisfy in the viscosity sense%
\begin{equation}
\left\vert
\begin{array}{l}
\dfrac{\partial \bar{u}}{\partial t}\left( t,x\right) +\tilde{f}\big (t,x,%
\bar{u}\left( t,x\right) ,D\bar{u}\left( t,x\right) ,D^{2}\bar{u}\left(
t,x\right) \big )\vspace*{2mm} \\
\text{\ \ \ \ \ \ }\;\text{\ \ \ \ \ \ \ }\;\text{\ \ \ \ }+e^{\lambda
t}\varphi _{-}^{\prime }\big (e^{-\lambda t}\bar{u}(t,x)\big )\leq 0\ \;%
\text{if}\;x\in \mathcal{D}\bigskip \\
\min \bigg\{\dfrac{\partial \bar{u}}{\partial t}\left( t,x\right) +\tilde{f}%
\big (t,x,\bar{u}\left( t,x\right) ,D\bar{u}\left( t,x\right) ,D^{2}\bar{u}%
\left( t,x\right) \big )\medskip \\
\text{\ \ \ \ }\left. +e^{\lambda t}\varphi _{-}^{\prime }\big (e^{-\lambda
t}\bar{u}(t,x)\big )\;,\;\big\langle\nabla \tilde{\ell}\left( x\right) ,D%
\bar{u}\left( t,x\right) \big\rangle+\delta \right. \vspace*{2mm} \\
-\tilde{g}\big (t,x,\bar{u}(t,x)\big )+e^{\lambda t}\psi _{-}^{\prime }\big (%
e^{-\lambda t}\bar{u}(t,x)\big )\bigg\}\leq 0\;\ \text{if}\;x\in Bd\left(
\mathcal{D}\right) .%
\end{array}%
\right.  \label{inegu}
\end{equation}%
Analogously we see that $\bar{v}$ satisfy in the viscosity sense:\medskip
\begin{equation}
\left\vert
\begin{array}{l}
\dfrac{\partial \bar{v}}{\partial t}\left( t,x\right) +\tilde{f}\big (t,x,%
\bar{v}\left( t,x\right) ,D\bar{v}\left( t,x\right) ,D^{2}\bar{v}\left(
t,x\right) \big )\vspace*{2mm} \\
\text{\ \ \ \ \ \ }\;\text{\ \ \ \ \ \ \ }\;\text{\ \ }+e^{\lambda t}\varphi
_{+}^{\prime }\big (e^{-\lambda t}\bar{v}(t,x)\big )-\varepsilon /t^{2}\geq
0\ \;\text{if}\;x\in \mathcal{D}\bigskip \\
\max \bigg\{\dfrac{\partial \bar{v}}{\partial t}\left( t,x\right) +\tilde{f}%
\big (t,x,\bar{v}\left( t,x\right) ,D\bar{v}\left( t,x\right) ,D^{2}\bar{v}%
\left( t,x\right) \big )\medskip \\
\left. +e^{\lambda t}\varphi _{-}^{\prime }\big (e^{-\lambda t}\bar{v}(t,x)%
\big )-\varepsilon /t^{2}\;,\;\big\langle\nabla \tilde{\ell}\left( x\right)
,D\bar{v}\left( t,x\right) \big\rangle-\delta \right. \vspace*{2mm} \\
-\tilde{g}\big (t,x,\bar{v}(t,x)\big )+e^{\lambda t}\psi _{+}^{\prime }\big (%
e^{-\lambda t}\bar{v}(t,x)\big )\bigg\}\geq 0\;\;\text{if}\;x\in Bd\left(
\mathcal{D}\right) ,%
\end{array}%
\right.  \label{inegv}
\end{equation}%
For simplicity of notation we continue to write $u,v$ for $\bar{u},\bar{v}$
respectively.

We assume now, to the contrary, that%
\begin{equation}
\underset{\left[ 0,T\right] \times \overline{\mathcal{D}}}{\max }\left(
u-v\right) ^{+}>0.  \label{ppra}
\end{equation}%
Exactly as in Theorem 4.2 in \cite{pa-ra/98} we have $(\hat{t},\hat{x})\in %
\left[ 0,T\right] \times Bd\left( \mathcal{D}\right) $ where $(\hat{t},\hat{x%
})$ is the maximum point, i.e.%
\begin{equation*}
u(\hat{t},\hat{x})-v(\hat{t},\hat{x})=\underset{\left[ 0,T\right] \times
\overline{\mathcal{D}}}{\max }\left( u-v\right) ^{+}>0.
\end{equation*}%
We put now (see also the proof of the Theorem 7.5 in Crandall, Ishii, Lions
\cite{cr-li/92})%
\begin{equation*}
\Phi _{n}\left( t,x,y\right) =u\left( t,x\right) -v\left( t,y\right) -\rho
_{n}\left( t,x,y\right) \text{, with }\left( t,x,y\right) \in \left[ 0,T%
\right] \times \overline{\mathcal{D}}\times \overline{\mathcal{D}},
\end{equation*}%
where%
\begin{equation}
\begin{array}{l}
\rho _{n}\left( t,x,y\right) =\dfrac{n}{2}\left\vert x-y\right\vert ^{2}+%
\tilde{g}\big (\hat{t},\hat{x},u(\hat{t},\hat{x})\big )\big\langle\nabla
\tilde{\ell}\left( \hat{x}\right) ,x-y\big\rangle+\left\vert x-\hat{x}%
\right\vert ^{4}\vspace*{2mm} \\
\;\;\;\;\;\;\;\;\;\;\;\;\;\;\;\;\;\;\;\;+|t-\hat{t}|^{4}-e^{\lambda \hat{t}%
}\psi _{-}^{\prime }\big (e^{-\lambda \hat{t}}u(\hat{t},\hat{x})\big )%
\big\langle\nabla \tilde{\ell}\left( \hat{x}\right) ,x-y\big\rangle.%
\end{array}
\label{defr}
\end{equation}%
Let it be $\left( t_{n},x_{n},y_{n}\right) $ a maximum point of $\Phi _{n}$.

We observe that $u\left( t,x\right) -v\left( t,x\right) -\left\vert x-\hat{x}%
\right\vert ^{4}-|t-\hat{t}|^{4}$ has in $(\hat{t},\hat{x})$ a unique
maximum point. Then, by Proposition 3.7 in Crandall, Ishii, Lions \cite%
{cr-li/92}, we have, as $n\rightarrow \infty $%
\begin{equation}
\begin{array}{l}
t_{n}\rightarrow \hat{t},\;x_{n}\rightarrow \hat{x},\;y_{n}\rightarrow \hat{x%
},\;n\left\vert x_{n}-y_{n}\right\vert ^{2}\rightarrow 0,\vspace*{2mm} \\
u\left( t_{n},x_{n}\right) \rightarrow u(\hat{t},\hat{x}),\;v\left(
t_{n},x_{n}\right) \rightarrow v(\hat{t},\hat{x}).%
\end{array}
\label{lpm}
\end{equation}%
But domain $\mathcal{D}$ verify the uniform exterior sphere condition:%
\begin{equation*}
\exists ~r_{0}>0\text{ such that }S\big (x+r_{0}\nabla \tilde{\ell}\left(
x\right) ,r_{0}\big )\cap \mathcal{D}=\emptyset \;,\;\text{ for }x\in
Bd\left( \mathcal{D}\right)
\end{equation*}%
where $S\left( x,r_{0}\right) $ denotes the closed ball of radius $r_{0}$
centered at $x$.

Then%
\begin{equation*}
\big |y-x-r_{0}\nabla \tilde{\ell}\left( x\right) \big |^{2}>r_{0}^{2}\text{
,\ \ for }x\in Bd\left( \mathcal{D}\right) \text{, }y\in \overline{\mathcal{D%
}},
\end{equation*}%
or equivalent%
\begin{equation}
\big\langle\nabla \tilde{\ell}\left( x\right) ,y-x\big\rangle<\frac{1}{2r_{0}%
}\left\vert y-x\right\vert ^{2}\text{ for }x\in Bd\left( \mathcal{D}\right)
\text{, }y\in \overline{\mathcal{D}}.  \label{se}
\end{equation}%
If we denote%
\begin{equation*}
B\left( t,x,r,q\right) =\big\langle\nabla \tilde{\ell}\left( x\right) ,q%
\big\rangle-\tilde{g}(t,x,r)
\end{equation*}%
then, if $x_{n}\in Bd\left( \mathcal{D}\right) $, we have, using the form of
$\rho _{n}$ given by (\ref{defr}) and (\ref{se}), that\newline
$%
\begin{array}{c}
B\big (t_{n},x_{n},u\left( t_{n},x_{n}\right) ,D_{x}\rho _{n}\left(
t_{n},x_{n},y_{n}\right) \big )=B\Big(t_{n},x_{n},u\left( t_{n},x_{n}\right)
,n\left( x_{n}-y_{n}\right) \vspace*{2mm}%
\end{array}%
$\newline
$%
\begin{array}{c}
+\tilde{g}\big (\hat{t},\hat{x},u(\hat{t},\hat{x})\big )\nabla \tilde{\ell}%
\left( \hat{x}\right) +4\left\vert x_{n}-\hat{x}\right\vert ^{2}\left( x_{n}-%
\hat{x}\right) -e^{\lambda \hat{t}}\psi _{-}^{\prime }\big (e^{-\lambda \hat{%
t}}u(\hat{t},\hat{x})\big )\nabla \tilde{\ell}\left( \hat{x}\right) \Big)%
\vspace*{2mm}%
\end{array}%
$\newline
$%
\begin{array}{c}
\geq -\dfrac{n}{2r_{0}}\left\vert x_{n}-y_{n}\right\vert ^{2}+\tilde{g}\big (%
\hat{t},\hat{x},u(\hat{t},\hat{x})\big )\big\langle\nabla \tilde{\ell}\left(
\hat{x}\right) ,\nabla \tilde{\ell}\left( x_{n}\right) \big\rangle\vspace*{%
2mm}%
\end{array}%
$\newline
$%
\begin{array}{c}
-\tilde{g}\big (t_{n},x_{n},u\left( t_{n},x_{n}\right) \big )+4\left\vert
x_{n}-\hat{x}\right\vert ^{2}\big\langle\nabla \tilde{\ell}\left(
x_{n}\right) ,x_{n}-\hat{x}\big\rangle\vspace*{2mm}%
\end{array}%
$\newline
$%
\begin{array}{c}
-e^{\lambda \hat{t}}\psi _{-}^{\prime }\big (e^{-\lambda \hat{t}}u(\hat{t},%
\hat{x})\big )\big\langle\nabla \tilde{\ell}\left( \hat{x}\right) ,\nabla
\tilde{\ell}\left( x_{n}\right) \big\rangle\vspace*{2mm}%
\end{array}%
$\newline
Then (\ref{lpm}) implies for $x_{n}\in Bd\left( \mathcal{D}\right) :$%
\begin{equation*}
\begin{array}{l}
\liminf\limits_{n\rightarrow \infty }\Big[B\left( t_{n},x_{n},u\left(
t_{n},x_{n}\right) ,D_{x}\rho _{n}\left( t_{n},x_{n},y_{n}\right) \right)
+\delta \vspace*{2mm} \\
+e^{\lambda t_{n}}\psi _{-}^{\prime }\big(e^{-\lambda t_{n}}u(t_{n},x_{n})%
\big)\Big]>0%
\end{array}%
\end{equation*}%
Analogously if $y_{n}\in Bd\left( \mathcal{D}\right) $ we infer%
\begin{equation*}
\begin{array}{l}
\limsup\limits_{n\rightarrow \infty }\Big[B\left( t_{n},y_{n},v\left(
t_{n},y_{n}\right) ,-D_{y}\rho _{n}\left( t_{n},x_{n},y_{n}\right) \right)
-\delta \vspace*{2mm} \\
+e^{\lambda t_{n}}\psi _{+}^{\prime }\big(e^{-\lambda t_{n}}v(t_{n},x_{n})%
\big)\Big]<0\text{.}%
\end{array}%
\end{equation*}%
Then from (\ref{inegu}), (\ref{inegv}) we conclude that%
\begin{equation}
\begin{array}{l}
p+\tilde{f}\big (t_{n},x_{n},u\left( t_{n},x_{n}\right) ,D_{x}\rho
_{n}\left( t_{n},x_{n},y_{n}\right) ,X\big )+e^{\lambda t_{n}}\varphi
_{-}^{\prime }\big (e^{-\lambda t_{n}}u(t_{n},x_{n})\big )\leq 0,\vspace*{2mm%
} \\
\text{\ \ \ \ \ \ }\;\text{\ \ \ \ \ \ \ }\;\text{\ \ \ \ \ \ \ }\;\text{\ \
\ \ \ \ \ }\;\text{\ for\ }\left( p,D_{x}\rho _{n}\left(
t_{n},x_{n},y_{n}\right) ,X\right) \in \overline{\mathcal{P}}%
^{2,+}u(t_{n},x_{n})%
\end{array}
\label{fneg}
\end{equation}%
and%
\begin{equation}
\begin{array}{l}
p+\tilde{f}\big (t_{n},y_{n},v\left( t_{n},y_{n}\right) ,-D_{y}\rho
_{n}\left( t_{n},x_{n},y_{n}\right) ,Y\big )\vspace*{2mm} \\
\;\;\;\;\;\;\;\;\;\;\;\;\;\;\;\;\;\;\;\;\;\;\;\;\;\;\;\;\;\;\;\;\;\;\;\;\;\;%
\;\;+e^{\lambda t_{n}}\varphi _{+}^{\prime }\big (e^{-\lambda
t_{n}}v(t_{n},y_{n})\big )\geq \dfrac{\varepsilon }{t^{2}}\;,\medskip \\
\text{\ \ \ \ }\;\;\;\;\;\;\;\;\;\;\;\;\;\;\;\;\;\;\;\;\;\;\;\;\;\;\;\;\;%
\text{for\ }\big (p,-D_{y}\rho _{n}\left( t_{n},x_{n},y_{n}\right) ,Y\big )%
\in \overline{\mathcal{P}}^{2,-}v(t_{n},y_{n})%
\end{array}
\label{fpoz}
\end{equation}%
From Theorem 8.3 in Crandall, Ishii, Lions \cite{cr-li/92} (apply with $k=2$%
, $\mathcal{O}_{1}=\mathcal{O}_{2}=\overline{\mathcal{D}}$, $u_{1}=u$, $%
u_{2}=-v$, $b_{1}=p$, $b_{2}=-p$ ) we deduce that there exists%
\begin{equation*}
\left( p,X,Y\right) \in \mathbb{R}\times \text{S\negthinspace }\mathbb{R}%
^{d\times d}\times \text{S\negthinspace }\mathbb{R}^{d\times d},
\end{equation*}%
such that%
\begin{equation*}
\begin{array}{c}
\big (p,D_{x}\rho _{n}\left( t_{n},x_{n},y_{n}\right) ,X\big )\in \overline{%
\mathcal{P}}^{2,+}u(t_{n},x_{n})\vspace*{2mm} \\
\big (p,-D_{y}\rho _{n}\left( t_{n},x_{n},y_{n}\right) ,Y\big )\in \overline{%
\mathcal{P}}^{2,-}v(t_{n},y_{n})%
\end{array}%
\end{equation*}%
and%
\begin{equation}
-\left( n+\left\vert \left\vert A\right\vert \right\vert \right)
\begin{pmatrix}
I & \ \ 0 \\
0 & \ \ I%
\end{pmatrix}%
\leq
\begin{pmatrix}
X & \ \ 0 \\
0 & -Y%
\end{pmatrix}%
\leq A+\frac{1}{n}A^{2},  \label{inegm}
\end{equation}%
where $A=D_{x,y}^{2}\rho _{n}\left( t_{n},x_{n},y_{n}\right) $. From (\ref%
{defr}) we have%
\begin{equation*}
\begin{array}{l}
A=n%
\begin{pmatrix}
\ \ I & \ -I \\
-I & \ \ \ I%
\end{pmatrix}%
+O\big (\left\vert x_{n}-\hat{x}\right\vert ^{2}\big ),\vspace*{2mm} \\
A^{2}=2n^{2}%
\begin{pmatrix}
\ \ I & \ -I \\
-I & \ \ \ I%
\end{pmatrix}%
+O\big (n\left\vert x_{n}-\hat{x}\right\vert ^{2}+\left\vert x_{n}-\hat{x}%
\right\vert ^{4}\big ),%
\end{array}%
\end{equation*}%
where $\left\vert O\left( h\right) \right\vert \leq C\left\vert h\right\vert
$ (the Landau symbol). Then (\ref{inegm}) become%
\begin{equation}
-\left( 3n+\kappa _{n}\right)
\begin{pmatrix}
I & \ \ 0 \\
0 & \ \ I%
\end{pmatrix}%
\leq
\begin{pmatrix}
X & \ \ 0 \\
0 & -Y%
\end{pmatrix}%
\leq 3n%
\begin{pmatrix}
\ \ I & \ -I \\
-I & \ \ \ I%
\end{pmatrix}%
+\kappa _{n}%
\begin{pmatrix}
I & \ \ 0 \\
0 & \ \ I%
\end{pmatrix}
\label{inegmat2}
\end{equation}%
where $\kappa _{n}\rightarrow 0.$ Now from (\ref{fneg}) and(\ref{fpoz})%
\begin{equation*}
\begin{array}{lll}
\dfrac{\varepsilon }{t^{2}} & \leq & \tilde{f}\big (t_{n},y_{n},v\left(
t_{n},y_{n}\right) ,-D_{y}\rho _{n}\left( t_{n},x_{n},y_{n}\right) ,Y\big )%
+e^{\lambda t_{n}}\varphi _{+}^{\prime }\big (e^{-\lambda
t_{n}}v(t_{n},y_{n})\big )\vspace*{2mm} \\
&  & -\tilde{f}\big (t_{n},x_{n},u\left( t_{n},x_{n}\right) ,D_{x}\rho
_{n}\left( t_{n},x_{n},y_{n}\right) ,X\big )-e^{\lambda t_{n}}\varphi
_{-}^{\prime }\big (e^{-\lambda t_{n}}u(t_{n},x_{n})\big )%
\end{array}%
\end{equation*}%
By (\ref{ppra}) and (\ref{lpm}) there exists $N\geq 1$ such that%
\begin{equation*}
u(t_{n},x_{n})>v(t_{n},y_{n}),\;\forall ~n\geq N
\end{equation*}%
and consequently%
\begin{equation*}
e^{\lambda t_{n}}\varphi _{-}^{\prime }\big (e^{-\lambda t}u(t_{n},x_{n})%
\big )\geq e^{\lambda t_{n}}\varphi _{+}^{\prime }\big (e^{-\lambda
t}v(t_{n},y_{n})\big )
\end{equation*}%
and%
\begin{equation*}
\begin{array}{l}
\tilde{f}\big (t_{n},y_{n},u\left( t_{n},x_{n}\right) ,-D_{y}\rho _{n}\left(
t_{n},x_{n},y_{n}\right) ,Y\big )\geq \medskip \\
\;\;\;\;\;\;\;\;\;\;\text{\ \ \ \ \ \ }\;\text{\ \ \ \ \ \ \ \ \ \ \ \ \ \ }%
\ \ \ \ \ \ \;\tilde{f}\big (t_{n},y_{n},v\left( t_{n},y_{n}\right)
,-D_{y}\rho _{n}\left( t_{n},x_{n},y_{n}\right) ,Y\big )%
\end{array}%
\end{equation*}%
Then, by definition (\ref{defft}) of $\tilde{f}$ and assumption (\ref{h7}),
we have
\begin{equation*}
\begin{array}{l}
\dfrac{\varepsilon }{\hat{t}^{2}}\leq \liminf_{n\rightarrow +\infty }\Big[%
\tilde{f}\big(t_{n},y_{n},u\left( t_{n},x_{n}\right) ,-D_{y}\rho _{n}\left(
t_{n},x_{n},y_{n}\right) ,Y\big)\medskip \\
\;\;\;\;\;\;-\tilde{f}\big(t_{n},x_{n},u\left( t_{n},x_{n}\right) ,D_{x}\rho
_{n}\left( t_{n},x_{n},y_{n}\right) ,X\big)\Big]\medskip \\
\;\;\;\;\;\leq \dfrac{1}{2}\mathrm{Tr}\big[\left( \sigma \sigma ^{\ast
}\right) (t_{n},x_{n})X-\left( \sigma \sigma ^{\ast }\right) (t_{n},y_{n})Y%
\big]%
\end{array}%
\end{equation*}%
But from (\ref{inegmat2}), $\forall $ $q,\tilde{q}\in \mathbb{R}^{d}$,%
\begin{equation*}
\left\langle Xq,q\right\rangle -\left\langle Y\tilde{q},\tilde{q}%
\right\rangle \leq 3n\left\vert q-\tilde{q}\right\vert ^{2}+\big (\left\vert
q\right\vert ^{2}+\left\vert \tilde{q}\right\vert ^{2}\big )\kappa _{n}.
\end{equation*}%
Hence%
\begin{align*}
& \mathrm{Tr}\big[\left( \sigma \sigma ^{\ast }\right) (t_{n},x_{n})X-\left(
\sigma \sigma ^{\ast }\right) (t_{n},y_{n})Y\big] \\
& =\sum\limits_{i=1}^{d}\big (X\sigma (t_{n},x_{n})e_{i},\sigma
(t_{n},x_{n})e_{i}\big )-\big (Y\sigma (t_{n},y_{n})e_{i},\sigma
(t_{n},y_{n})e_{i}\big ) \\
& \leq 3C~n\left\vert x_{n}-y_{n}\right\vert ^{2}+\big (\left\vert \sigma
(t_{n},x_{n})\right\vert ^{2}+\left\vert \sigma (t_{n},y_{n})\right\vert ^{2}%
\big )\kappa _{n},
\end{align*}%
and consequently%
\begin{equation*}
\dfrac{\varepsilon }{\hat{t}^{2}}\leq 0
\end{equation*}%
that is a contradiction.

Then%
\begin{equation*}
u\left( t,x\right) \leq v\left( t,x\right) ,\;\forall \left( t,x\right) \in
\left[ 0,T\right] \times \overline{\mathcal{D}}.
\end{equation*}
\end{proof}

\section{Backward stochastic variational inequalities}

Let $\left\{ W_{t}:t\geq0\right\} $ be a $d$-dimensional standard Brownian
motion defined on some complete probability space $(\Omega,\mathcal{F},%
\mathbb{P})$. We denote by $\left\{ \mathcal{F}_{t}:t\geq0\right\} $ the
natural filtration generated by $\left\{ W_{t}:t\geq0\right\} $ \ and
augmented by $\mathcal{N}$ the set of$\;\mathbb{P}$- null events of $%
\mathcal{F}$:
\begin{equation*}
\mathcal{F}_{t}=\sigma\{W_{r}:0\leq r\leq t\}\vee\mathcal{N}.
\end{equation*}

Let $\tau:\Omega\rightarrow\left[ 0,\infty\right) \;$be an a.s. $\mathcal{F}%
_{t}$-stopping time and $\smallskip$

\begin{itemize}
\item $\left\{ A_{t}:t\geq0\right\} $ be a continuous one-dimensional
increasing progressively measurable stochastic process (p.m.s.p.) satisfying
$A_{0}=0$.\medskip
\end{itemize}

We shall study the existence and uniqueness of a solution $\left( Y,Z\right)
$ of the following backward stochastic variational inequality (BSVI)~:%
\begin{equation}
\left\{
\begin{array}{r}
dY_{t}+F\left( t,Y_{t},Z_{t}\right) dt+G\left( t,Y_{t}\right)
dA_{t}\in\partial\varphi\left( Y_{t}\right) dt+\partial\psi\left(
Y_{t}\right) dA_{t}+Z_{t}dW_{t},\medskip \\
0\leq t\leq\tau, \\
\multicolumn{1}{l}{Y_{\tau}=\xi.}%
\end{array}
\right.  \label{bsvi}
\end{equation}

\subsection{Assumptions and results}

Let $\lambda,\mu\geq0$.

Let%
\begin{equation*}
\mathcal{H}_{k}^{\lambda,\mu}\subset L^{2}(\mathbb{R}_{+}\times\Omega
,e^{^{\lambda s+\mu A_{s}}}\mathbf{1}_{\left[ 0,\tau\right] }\left( s\right)
ds\otimes d\mathbb{P};\mathbb{R}^{k})
\end{equation*}
the Hilbert space of p.m.s.p. $f:\Omega\times\left[ 0,\infty\right)
\rightarrow\mathbb{R}^{k}$ such that%
\begin{equation*}
\left\Vert f\right\Vert _{\mathcal{H}}=\bigg [\mathbb{E}\Big (\int_{0}^{\tau
}\text{ }e^{^{\lambda s+\mu A_{s}}}\big |f\left( s\right) \big |^{2}ds\Big )%
\bigg ]^{1/2}<\infty
\end{equation*}
and%
\begin{equation*}
\mathcal{\tilde{H}}_{k}^{\lambda,\mu}\subset L^{2}(\mathbb{R}_{+}\times
\Omega,e^{^{\lambda s+\mu A_{s}}}\mathbf{1}_{\left[ 0,\tau\right] }\left(
s\right) dA_{s}\otimes d\mathbb{P};\mathbb{R}^{k})
\end{equation*}
the Hilbert space of p.m.s.p. $f:\Omega\times\left[ 0,\infty\right)
\rightarrow\mathbb{R}^{k}$ such that%
\begin{equation*}
\left\Vert f\right\Vert _{\mathcal{\tilde{H}}}=\bigg[\mathbb{E}\Big(\int
_{0}^{\tau}\text{ }e^{^{\lambda s+\mu A_{s}}}\big |f\left( s\right) \big |%
^{2}dA_{s}\Big)\bigg]^{1/2}<\infty\;.
\end{equation*}

We also introduce the notation $\mathcal{S}_{k}^{\lambda,\mu}$ for the
Banach space of p.m.s.p.\ $f:\Omega\times\left[ 0,\infty\right) \rightarrow
\mathbb{R}^{k}$ such that
\begin{equation*}
\left\Vert f\right\Vert _{\mathcal{S}}=\bigg [\mathbb{E}\Big (\underset{%
0\leq t\leq\tau}{\sup}\text{ }e^{\lambda t+\mu A_{t}}\big |f\left( t\right) %
\big |^{2}\Big )\bigg ]^{1/2}<\infty\;.
\end{equation*}

With respect to BSVI (\ref{bsvi}) we formulate the following assumptions:

\begin{itemize}
\item Let $F:\Omega\times\left[ 0,\infty\right) \times\mathbb{R}^{k}\times%
\mathbb{R}^{k\times d}\rightarrow\mathbb{R}^{k}$, $G:\Omega\times\left[
0,\infty\right) \times\mathbb{R}^{k}\rightarrow\mathbb{R}^{k}$ satisfy that
there exist $\alpha,\beta\in\mathbb{R}$, $L\geq0$ and $\eta,\gamma
:[0,\infty)\times\Omega\rightarrow\lbrack0,\infty)$ an p.m.s.p. such that
for all $t\geq0$, $y,y^{\prime}\in\mathbb{R}^{k}$, $z,z^{\prime}\in\mathbb{R}%
^{k\times d}:$%
\begin{equation}
\left\vert \;\;%
\begin{array}{cl}
\left( i\right) & F\left( \cdot,\cdot,y,z\right) \;\text{is p.m.s.p.,}%
\medskip \\
\left( ii\right) & y\longrightarrow F\left( \omega,t,y,z\right) :\mathbb{R}%
^{k}\rightarrow\mathbb{R}^{k}~\text{is continuous,}\;\text{a.s.}\medskip \\
\left( iii\right) & \big\langle y-y^{\prime},F\left( t,y,z\right) -F\left(
t,y^{\prime},z\right) \big\rangle\leq\alpha\left\vert y-y^{\prime
}\right\vert ^{2},\;\text{a.s.}\medskip \\
\left( iv\right) & \big|F\left( t,y,z\right) -F\left( t,y,z^{\prime }\right) %
\big|\leq L\left\Vert z-z^{\prime}\right\Vert ,\;\text{a.s.}\medskip \\
\left( v\right) & \left\vert F\left( t,y,z\right) \right\vert \leq\eta _{t}+L%
\big(\left\vert y\right\vert +\left\Vert z\right\Vert \big),\;\text{a.s.}%
\end{array}
\right.  \label{ipF}
\end{equation}
and%
\begin{equation}
\left\vert \;\;%
\begin{array}{cl}
\left( i\right) & G\left( \cdot,\cdot,y\right) \;\text{is p.m.s.p.,}\medskip
\\
\left( ii\right) & y\longrightarrow G\left( \omega,t,y\right) :\mathbb{R}%
^{k}\rightarrow\mathbb{R}^{k}\;\text{is continuous,}\;\text{a.s.}\medskip \\
\left( iii\right) & \big\langle y-y^{\prime},G\left( t,y\right) -G\left(
t,y^{\prime}\right) \big\rangle\leq\beta\left\vert y-y^{\prime}\right\vert
^{2},\;\text{a.s.}\medskip \\
\left( iv\right) & \big|G\left( t,y\right) \big|\leq\gamma_{t}+L\left\vert
y\right\vert ,\;\text{a.s.}%
\end{array}
\right.  \label{ipG}
\end{equation}

\item Terminal date $\xi $ is an $\mathbb{R}^{k}$-valued $\mathcal{F}_{\tau
} $-measurable random variable such that there exists $\lambda >2\alpha
+2L^{2}+1$, $\mu >2\beta +1:$%
\begin{equation}
\begin{array}{l}
M\left( \tau \right) \overset{def}{=}\mathbb{E}e^{\lambda \tau +\mu A_{\tau
}}\big(\left\vert \xi \right\vert ^{2}+\varphi \left( \xi \right) +\psi
\left( \xi \right) \big)\medskip \\
\displaystyle\quad \quad \quad \quad +\mathbb{E}\int_{0}^{\tau }e^{\lambda
s+\mu A_{s}}\big[\left\vert \eta _{s}\right\vert ^{2}ds+\left\vert \gamma
_{s}\right\vert ^{2}dA_{s}\big]<\infty .%
\end{array}
\label{ipxi}
\end{equation}

\item Let it be $\varphi,\psi$ such that
\begin{equation}
\left\vert \;\;%
\begin{array}{cl}
\left( i\right) & \varphi,\psi:\mathbb{R}^{k}\rightarrow(-\infty ,+\infty]%
\text{ \ are proper convex l.s.c. functions,}\medskip \\
\left( ii\right) & \varphi\left( y\right) \geq\varphi\left( 0\right)
=0,\;\psi\left( y\right) \geq\psi\left( 0\right) =0,%
\end{array}
\right.  \label{ipfi}
\end{equation}
The subdifferentials are defined by%
\begin{equation*}
\partial\varphi\left( x\right) =\left\{ v\in\mathbb{R}^{k}:\left\langle
v,y-x\right\rangle +\varphi\left( x\right) \leq\varphi\left( y\right)
,\;\forall y\in\mathbb{R}^{k}\right\}
\end{equation*}
and similar for $\psi.$
\end{itemize}

The existence result for (\ref{bsvi}) will be obtained via Yosida
approximations. Define for $\varepsilon>0$ the convex $C^{1}$-function $%
\varphi_{\varepsilon}$ by%
\begin{equation*}
\varphi_{\varepsilon}\left( y\right) =\inf\left\{ \frac{1}{2\varepsilon }%
\left\vert y-v\right\vert ^{2}+\varphi\left( v\right) :v\in\mathbb{R}%
^{k}\right\}
\end{equation*}
(and similar for $\psi_{\varepsilon}$).

Denoting%
\begin{equation*}
J_{\varepsilon }y=\left( I+\varepsilon \partial \varphi \right) ^{-1}\left(
y\right) \;\;\text{and \ \ }\nabla \varphi _{\varepsilon }\left( y\right) =%
\frac{y-J_{\varepsilon }y}{\varepsilon }~.
\end{equation*}%
Hence $y\rightarrow \nabla \varphi _{\varepsilon }\left( y\right) $ is an
monotone Lipschitz function and%
\begin{equation*}
\varphi _{\varepsilon }\left( y\right) =\frac{1}{2\varepsilon }\left\vert
y-J_{\varepsilon }y\right\vert ^{2}+\varphi \left( J_{\varepsilon }y\right) .
\end{equation*}%
(analog for $\psi _{\varepsilon }$).

\begin{itemize}
\item We introduce now the \textit{compatibility assumptions }:\newline
for all \ $\varepsilon>0$, $t\geq0$, $y\in\mathbb{R}^{k}$ and $z\in\mathbb{R}%
^{k\times d}$%
\begin{equation}
\left\vert \;\;%
\begin{array}{rl}
\left( i\right) \; & \big\langle\nabla\varphi_{\varepsilon}\left( y\right)
,\nabla\psi_{\varepsilon}\left( y\right) \big\rangle\geq0,\medskip \\
\left( ii\right) \; & \big\langle\nabla\varphi_{\varepsilon}\left( y\right)
,G\left( t,y\right) \big\rangle\leq\big\langle\nabla
\psi_{\varepsilon}\left( y\right) ,G\left( t,y\right) \big\rangle%
^{+},\medskip \\
\left( iii\right) \; & \big\langle\nabla\psi_{\varepsilon}\left( y\right)
,F\left( t,y,z\right) \big\rangle\leq\big\langle\nabla\varphi_{\varepsilon
}\left( y\right) ,F\left( t,y,z\right) \big\rangle^{+}.%
\end{array}
\right.  \label{comp}
\end{equation}
\end{itemize}

\begin{definition}
$\left( Y,Z,U,V\right) $ will be called a solution of BSVI (\ref{bsvi}) if%
\begin{equation}
\begin{array}{ll}
\left( a\right) \; & \displaystyle Y\in\mathcal{S}_{k}^{\lambda,\mu}\cap%
\mathcal{H}_{k}^{\lambda,\mu}\cap\mathcal{\tilde{H}}_{k}^{\lambda,\mu
},\;Z\in\mathcal{H}_{k\times d}^{\lambda,\mu},\medskip \\
\left( b\right) \; & U\in\mathcal{H}_{k}^{\lambda,\mu},\;V\in\mathcal{%
\tilde {H}}_{k}^{\lambda,\mu},\medskip \\
\left( c\right) \; & \displaystyle\mathbb{E}\int_{0}^{\tau}e^{\lambda s+\mu
A_{s}}\big(\varphi\left( Y_{s}\right) ds+\psi\left( Y_{s}\right) dA_{s}\big)%
<\infty,\medskip \\
\left( d\right) \; & \left( Y_{t},U_{t}\right) \in\partial\varphi ,\;\mathbb{%
P}\left( d\omega\right) \otimes dt\;,\;\;\left( Y_{t},V_{t}\right)
\in\partial\psi,\;\mathbb{P}\left( d\omega\right) \otimes A\left(
\omega,dt\right) \\
& \;\;\;\text{\ \ \ }\;\;\;\text{\ \ \ }\;\;\;\text{\ \ \ }\;\;\;\text{\ \ \
}\;\;\;\text{\ \ \ }\;\;\;\text{\ \ \ }\;\;\;\text{\ \ \ }\;\;\;\text{\ \ \ }%
\;\;\;\text{\ \ \ }\;\;\;\text{\ \ \ }\;\;\;\text{\ \ \ }\;\text{\ a.e. on }%
\Omega\times\left[ 0,\tau\right] ,\medskip \\
\left( e\right) \; & \displaystyle Y_{t}+\int_{t\wedge\tau}^{\tau}U_{s}ds+%
\int_{t\wedge\tau}^{\tau}V_{s}dA_{s}=\xi+\int_{t\wedge\tau}^{\tau }F\left(
s,Y_{s},Z_{s}\right) ds\medskip \\
& \;\;\;\;\;\;\;\;\;\ \ \ \ \ \ \ \ \ \ \;\;\;\displaystyle%
+\int_{t\wedge\tau }^{\tau}G\left( s,Y_{s}\right)
dA_{s}-\int_{t\wedge\tau}^{\tau}Z_{s}dW_{s},\;\;\text{for all }t\geq0\;\text{%
a.s.}\medskip%
\end{array}
\label{def}
\end{equation}
\end{definition}

In all that follows, $C$ denotes a constant, which may depend only on $%
\mu,\alpha,\beta$ and $L$, which may vary from line to line.

\begin{proposition}
\label{p1}Let assumptions (\ref{ipF}), (\ref{ipG}), and (\ref{ipfi}). If $%
\left( Y,Z,U,V\right) $ and $(\tilde{Y},\tilde{Z},\tilde{U},\tilde{V})$ are
corresponding solutions to $\xi $ and $\tilde{\xi}$ which satisfy (\ref{ipxi}%
), then%
\begin{equation}
\begin{array}{r}
\displaystyle\mathbb{E}\int_{0}^{\tau }e^{\lambda s+\mu A_{s}}\left[ |Y_{s}-%
\tilde{Y}_{s}|^{2}\left( ds+dA_{s}\right) +||Z_{s}-\tilde{Z}_{s}||^{2}ds%
\right] \medskip \\
+\mathbb{E}\underset{0\leq t\leq \tau }{\sup }e^{\lambda t+\mu A_{t}}|Y_{t}-%
\tilde{Y}_{t}|^{2}\leq C~\mathbb{E}\left[ e^{\lambda \tau +\mu A_{\tau
}}|\xi -\tilde{\xi}|^{2}\right] .%
\end{array}
\label{unic}
\end{equation}
\end{proposition}

\begin{proof}
From It\^{o}'s formula we have\newline
$%
\begin{array}{c}
\displaystyle e^{\lambda \left( t\wedge \tau \right) +\mu A_{t\wedge \tau
}}|Y_{t\wedge \tau }-\tilde{Y}_{t\wedge \tau }|^{2}+\int_{t\wedge \tau
}^{\tau }e^{\lambda s+\mu A_{s}}|Y_{s}-\tilde{Y}_{s}|^{2}\left( \lambda
ds+\mu dA_{s}\right) \medskip%
\end{array}%
$\newline
$%
\begin{array}{c}
\displaystyle+2\int_{t\wedge \tau }^{\tau }e^{\lambda s+\mu A_{s}}\langle
Y_{s}-\tilde{Y}_{s},U_{s}-\tilde{U}_{s}\rangle ds+2\int_{t\wedge \tau
}^{\tau }e^{\lambda s+\mu A_{s}}\langle Y_{s}-\tilde{Y}_{s},V_{s}-\tilde{V}%
_{s}\rangle dA_{s}\medskip%
\end{array}%
$\newline
$%
\begin{array}{c}
\displaystyle+\int_{t\wedge \tau }^{\tau }e^{\lambda s+\mu A_{s}}||Z_{s}-%
\tilde{Z}_{s}||^{2}ds\medskip%
\end{array}%
$\newline
$%
\begin{array}{c}
\displaystyle=e^{\lambda \tau +\mu A_{\tau }}|\xi -\tilde{\xi}%
|^{2}+2\int_{t\wedge \tau }^{\tau }e^{\lambda s+\mu A_{s}}\big\langle Y_{s}-%
\tilde{Y}_{s},F\left( s,Y_{s},Z_{s}\right) -F(s,\tilde{Y}_{s},\tilde{Z}_{s})%
\big\rangle ds\medskip%
\end{array}%
$\newline
$%
\begin{array}{c}
\displaystyle+2\int_{t\wedge \tau }^{\tau }e^{\lambda s+\mu A_{s}}%
\big\langle
Y_{s}-\tilde{Y}_{s},G\left( s,Y_{s}\right) -G(s,\tilde{Y}_{s})\big\rangle %
dA_{s}\medskip%
\end{array}%
$\newline
$%
\begin{array}{c}
\displaystyle-2\int_{t\wedge \tau }^{\tau }e^{\lambda s+\mu A_{s}}%
\big\langle
Y_{s}-\tilde{Y}_{s},(Z_{s}-\tilde{Z}_{s})dW_{s}\big\rangle\medskip%
\end{array}%
$\newline
Since%
\begin{equation*}
\langle Y_{s}-\tilde{Y}_{s},U_{s}-\tilde{U}_{s}\rangle ds\geq 0,\;\quad
\langle Y_{s}-\tilde{Y}_{s},V_{s}-\tilde{V}_{s}\rangle dA_{s}\geq 0,
\end{equation*}%
\begin{equation*}
2\big\langle Y_{s}-\tilde{Y}_{s},F\left( s,Y_{s},Z_{s}\right) -F(s,\tilde{Y}%
_{s},\tilde{Z}_{s})\big\rangle\leq (2\alpha +2L^{2}+1)|Y_{s}-\tilde{Y}%
_{s}|^{2}+\dfrac{1}{2}||Z_{s}-\tilde{Z}_{s}||^{2}
\end{equation*}%
and%
\begin{equation*}
2\big\langle Y_{s}-\tilde{Y}_{s},G\left( s,Y_{s}\right) -G(s,\tilde{Y}_{s})%
\big\rangle\leq (2\beta +1)|Y_{s}-\tilde{Y}_{s}|^{2}
\end{equation*}%
then (using also the Burkholder--Davis--Gundy's inequality), inequality (\ref%
{unic}) follows.\hfill $\bigskip $
\end{proof}

The main result of this section is given by

\begin{theorem}
\label{t1}Let assumptions (\ref{ipF})-(\ref{comp}) be satisfied. Then there
exists a unique solution $\left( Y,Z,U,V\right) $ for (\ref{bsvi}).
\end{theorem}

\subsection{BSVI - proof of the existence}

Consider the approximating equation%
\begin{equation}
\begin{array}{l}
\displaystyle Y_{t}^{\varepsilon }+\int_{t\wedge \tau }^{\tau }\nabla
\varphi _{\varepsilon }\left( Y_{s}^{\varepsilon }\right) ds+\int_{t\wedge
\tau }^{\tau }\nabla \psi _{\varepsilon }\left( Y_{s}^{\varepsilon }\right)
dA_{s}=\xi +\int_{t\wedge \tau }^{\tau }F\left( s,Y_{s}^{\varepsilon
},Z_{s}^{\varepsilon }\right) ds\medskip \\
\;\;\ \ \ \ \ \ \ \ \ \ \ \ \ \ \ \ \ \ \ \ \ \ \ \ \ \ \ \ \ \ \ %
\displaystyle+\int_{t\wedge \tau }^{\tau }G\left( s,Y_{s}^{\varepsilon
}\right) dA_{s}-\int_{t\wedge \tau }^{\tau }Z_{s}^{\varepsilon
}dW_{s},\forall t\geq 0,P-a.s.%
\end{array}
\label{ecapr}
\end{equation}%
Since $\nabla \varphi _{\varepsilon },\nabla \psi _{\varepsilon }:\mathbb{R}%
^{k}\rightarrow \mathbb{R}^{k}$ are Lipschitz functions then, by a standard
argument (Banach fixed point theorem when $y\rightarrow F\left( t,y,z\right)
$ and $y\rightarrow G\left( t,y\right) $ are uniformly Lipschitz functions
and Lipschitz approximations when $y\rightarrow \alpha y-F\left(
t,y,z\right) $ and $y\rightarrow \beta y-G\left( t,y\right) $ are continuous
monotone functions) (see also \cite{pa-zh/98}), equation (\ref{ecapr}) has a
unique solution%
\begin{equation*}
\left( Y^{\varepsilon },Z^{\varepsilon }\right) \in \big(\mathcal{S}%
_{k}^{\lambda ,\mu }\cap \mathcal{H}_{k}^{\lambda ,\mu }\cap \mathcal{\tilde{%
H}}_{k}^{\lambda ,\mu }\big)\times \mathcal{H}_{k\times d}^{\mu ,\lambda }
\end{equation*}

\begin{proposition}
\label{p2}Let assumptions (\ref{ipF})-(\ref{ipfi}) be satisfied. Then%
\begin{equation}
\begin{array}{r}
\displaystyle\mathbb{E}\bigg[\underset{0\leq t\leq \tau }{\sup }e^{\lambda
t+\mu A_{t}}\left\vert Y_{t}^{\varepsilon }\right\vert ^{2}+\int_{0}^{\tau }%
\text{ }e^{^{\lambda s+\mu A_{s}}}\big(\left\vert Y_{s}^{\varepsilon
}\right\vert ^{2}+\left\Vert Z_{s}^{\varepsilon }\right\Vert ^{2}\big)%
ds\medskip \\
\displaystyle+\int_{0}^{\tau }\text{ }e^{^{^{\lambda s+\mu
A_{s}}}}\left\vert Y_{s}^{\varepsilon }\right\vert ^{2}dA_{s}\bigg]\leq
C~M\left( \tau \right)%
\end{array}
\label{ineg1}
\end{equation}
\end{proposition}

\begin{proof}
It\^{o}'s formula for $e^{\lambda t+\mu A_{t}}\left\vert Y_{t}^{\varepsilon
}\right\vert ^{2}$ yields\newline
$%
\begin{array}{c}
\displaystyle e^{\lambda \left( t\wedge \tau \right) +\mu A_{t\wedge \tau
}}\left\vert Y_{t\wedge \tau }^{\varepsilon }\right\vert ^{2}+\int_{t\wedge
\tau }^{\tau }e^{\lambda s+\mu A_{s}}\left\vert Y_{s}^{\varepsilon
}\right\vert ^{2}\left( \lambda ds+\mu dA_{s}\right) +\int_{t\wedge \tau
}^{\tau }e^{\lambda s+\mu A_{s}}\left\Vert Z_{s}^{\varepsilon }\right\Vert
^{2}ds\medskip%
\end{array}%
$\newline
$%
\begin{array}{c}
\displaystyle+2\int_{t\wedge \tau }^{\tau }e^{\lambda s+\mu A_{s}}\Big[%
\big\langle Y_{s}^{\varepsilon },\nabla \varphi _{\varepsilon }\left(
Y_{s}^{\varepsilon }\right) \big\rangle\lambda ds+\big\langle %
Y_{s}^{\varepsilon },\nabla \psi _{\varepsilon }\left( Y_{s}^{\varepsilon
}\right) \big\rangle\mu dA_{s}\Big]=e^{\lambda \tau +\mu A_{\tau
}}\left\vert \xi \right\vert ^{2}\medskip%
\end{array}%
$\newline
$%
\begin{array}{c}
\displaystyle+2\int_{t\wedge \tau }^{\tau }e^{\lambda s+\mu A_{s}}%
\big\langle
Y_{s}^{\varepsilon },F\left( s,Y_{s}^{\varepsilon },Z_{s}^{\varepsilon
}\right) \big\rangle ds+2\int_{t\wedge \tau }^{\tau }e^{\lambda s+\mu A_{s}}%
\big\langle Y_{s}^{\varepsilon },G\left( s,Y_{s}^{\varepsilon }\right) %
\big\rangle dA_{s}\medskip%
\end{array}%
$\newline
$%
\begin{array}{c}
\displaystyle-2\int_{t\wedge \tau }^{\tau }e^{\lambda s+\mu
A_{s}}\left\langle Y_{s}^{\varepsilon },Z_{s}^{\varepsilon
}dW_{s}\right\rangle .\medskip%
\end{array}%
$\newline
But from Schwartz's inequality and assumptions (\ref{ipF})-(\ref{ipxi}) we
obtain%
\begin{equation*}
\begin{array}{l}
2\big\langle Y_{s},F\left( s,Y_{s},Z_{s}\right) \big\rangle\leq 2\alpha
\left\vert Y_{s}\right\vert ^{2}+2L\left\vert Y_{s}\right\vert \left\Vert
Z_{s}\right\Vert +2\left\vert Y_{s}\right\vert \big|F\left( s,0,0\right) %
\big|\medskip \\
\;\;\;\;\;\;\;\;\;\;\;\;\;\;\;\;\;\ \;\ \ \ \ \;\;\;\
\;\;\;\;\;\;\;\;\;\;\;\leq \left( 2\alpha +2L^{2}+1\right) \left\vert
Y_{s}\right\vert ^{2}+\dfrac{1}{2}\left\Vert Z_{s}\right\Vert ^{2}+\big|%
F\left( s,0,0\right) \big|^{2}\medskip%
\end{array}%
\end{equation*}%
and%
\begin{equation*}
2\big\langle Y_{s},G\left( s,Y_{s},Z_{s}\right) \big\rangle\leq 2\beta
\left\vert Y_{s}\right\vert ^{2}+2\left\vert Y_{s}\right\vert \big|G\left(
s,0\right) \big|\leq \left( 2\beta +1\right) \left\vert Y_{s}\right\vert
^{2}+\big|G\left( s,0\right) \big|^{2}
\end{equation*}%
Hence, using also that $\left\langle y,\nabla \varphi _{\varepsilon }\left(
y\right) \right\rangle ,\left\langle y,\nabla \psi _{\varepsilon }\left(
y\right) \right\rangle \geq 0,\medskip $\newline
$%
\begin{array}{c}
\displaystyle e^{\lambda \left( t\wedge \tau \right) +\mu A_{t\wedge \tau
}}\left\vert Y_{t\wedge \tau }^{\varepsilon }\right\vert ^{2}+\int_{t\wedge
\tau }^{\tau }e^{\lambda s+\mu A_{s}}\left\vert Y_{s}^{\varepsilon
}\right\vert ^{2}(\lambda -2\alpha -2L^{2}-1)ds+\medskip%
\end{array}%
$\newline
$%
\begin{array}{c}
\displaystyle\int_{t\wedge \tau }^{\tau }e^{\lambda s+\mu A_{s}}\left\vert
Y_{s}^{\varepsilon }\right\vert ^{2}\left( \mu -2\beta -1\right) dA_{s}+%
\dfrac{1}{2}\int_{t\wedge \tau }^{\tau }e^{\lambda s+\mu A_{s}}\left\Vert
Z_{s}^{\varepsilon }\right\Vert ^{2}ds\leq e^{\lambda \tau +\mu A_{\tau
}}\left\vert \xi \right\vert ^{2}\medskip%
\end{array}%
$\newline
$%
\begin{array}{c}
\displaystyle+\int_{t\wedge \tau }^{\tau }e^{\lambda s+\mu A_{s}}\Big(\big|%
F\left( s,0,0\right) \big|^{2}ds+\big|G\left( s,0\right) \big|^{2}dA_{s}\Big)%
\medskip%
\end{array}%
$\newline
$%
\begin{array}{c}
\displaystyle-2\int_{t\wedge \tau }^{\tau }e^{\lambda s+\mu
A_{s}}\left\langle Y_{s}^{\varepsilon },Z_{s}^{\varepsilon
}dW_{s}\right\rangle \medskip%
\end{array}%
$\newline
that clearly yields (for $\lambda >2\alpha +2L^{2}+1$ and $\mu >2\beta +1$):%
\begin{equation*}
\mathbb{E}\int_{0}^{\tau }\text{ }e^{^{\lambda s+\mu A_{s}}}\big[\left\vert
Y_{s}^{\varepsilon }\right\vert ^{2}\left( ds+dA_{s}\right) +\left\Vert
Z_{s}^{\varepsilon }\right\Vert ^{2}ds\big]\leq C~M\left( \tau \right)
\end{equation*}%
Since, by Burkholder--Davis--Gundy's inequality,%
\begin{equation*}
\begin{array}{r}
\displaystyle\mathbb{E}\underset{0\leq t\leq \tau }{\sup }\left\vert
\int_{t\wedge \tau }^{\tau }e^{\lambda s+\mu A_{s}}\left\langle
Y_{s}^{\varepsilon },Z_{s}^{\varepsilon }dW_{s}\right\rangle \right\vert
\leq 3\mathbb{E}\left( \int_{0}^{\tau }e^{2\left( \lambda s+\mu A_{s}\right)
}\big|\left\langle Y_{s}^{\varepsilon },Z_{s}^{\varepsilon }\right\rangle %
\big|^{2}ds\right) ^{1/2} \\
\displaystyle\leq \dfrac{1}{4}\mathbb{E}\underset{0\leq t\leq \tau }{\sup }%
e^{\lambda t+\mu A_{t}}\left\vert Y_{t}^{\varepsilon }\right\vert ^{2}+C~%
\mathbb{E}\int_{0}^{\tau }\text{ }e^{^{\lambda s+\mu A_{s}}}\left\Vert
Z_{s}^{\varepsilon }\right\Vert ^{2}ds\;,%
\end{array}%
\end{equation*}%
then it follows
\begin{equation*}
\mathbb{E}\underset{0\leq t\leq \tau }{\sup }e^{\lambda t+\mu
A_{t}}\left\vert Y_{t}^{\varepsilon }\right\vert ^{2}\leq C~M\left( \tau
\right) ~.
\end{equation*}%
The proof is complete.\hfill
\end{proof}

\begin{proposition}
\label{p3}Let assumptions (\ref{ipF})-(\ref{comp}) be satisfied. Then there
exists a positive constant $C$ such that for any stopping time $\theta \in %
\left[ 0,\tau \right] :$%
\begin{equation}
\begin{array}{cl}
\left( a\right) & \displaystyle\mathbb{E}\int_{0}^{\tau }e^{\lambda s+\mu
A_{s}}\Big(\big|\nabla \varphi _{\varepsilon }\left( Y_{s}^{\varepsilon
}\right) \big|^{2}ds+\big|\nabla \psi _{\varepsilon }\left(
Y_{s}^{\varepsilon }\right) \big|^{2}dA_{s}\Big)\leq C~M\left( \tau \right)
\;,\medskip \\
\left( b\right) & \displaystyle\mathbb{E}\int_{0}^{\tau }e^{\lambda s+\mu
A_{s}}\Big(\varphi \big(J_{\varepsilon }\left( Y_{s}^{\varepsilon }\right) %
\big)ds+\psi \big(\hat{J}_{\varepsilon }\left( Y_{s}^{\varepsilon }\right) %
\big)dA_{s}\Big)\leq C~M\left( \tau \right) \;,\bigskip \\
\left( c\right) & \mathbb{E}e^{\lambda \theta +\mu A_{\theta }}\Big(\big|%
Y_{\theta }^{\varepsilon }-J_{\varepsilon }\left( Y_{\theta }^{\varepsilon
}\right) \big|^{2}+\big|Y_{\theta }^{\varepsilon }-\hat{J}_{\varepsilon
}\left( Y_{\theta }^{\varepsilon }\right) \big|^{2}\Big)\leq \varepsilon
~C~M\left( \tau \right) \;,\bigskip \\
\left( d\right) & \mathbb{E}e^{\lambda \theta +\mu A_{\theta }}\Big(\varphi %
\big(J_{\varepsilon }\left( Y_{\theta }^{\varepsilon }\right) \big)+\psi %
\big(\hat{J}_{\varepsilon }\left( Y_{\theta }^{\varepsilon }\right) \big)%
\Big)\leq C~M\left( \tau \right) \;.%
\end{array}
\label{ineg2}
\end{equation}
\end{proposition}

\begin{proof}
Essential for the proof is the stochastic subdifferential inequality
introduced by Pardoux and R\u{a}\c{s}canu in \cite{pa-ra/98}, 1998. We will
use this inequality for our purpose. First we write the subdifferential
inequality%
\begin{equation*}
\begin{array}{r}
e^{\lambda s+\mu A_{s}}\varphi _{\varepsilon }\left( Y_{s}^{\varepsilon
}\right) \geq (e^{\lambda s+\mu A_{s}}-e^{\lambda r+\mu A_{r}})\varphi
_{\varepsilon }\left( Y_{s}^{\varepsilon }\right) +e^{\lambda r+\mu
A_{r}}\varphi _{\varepsilon }\left( Y_{r}^{\varepsilon }\right) \medskip \\
+e^{\lambda r+\mu A_{r}}\big\langle\nabla \varphi _{\varepsilon }\left(
Y_{r}^{\varepsilon }\right) ,Y_{s}^{\varepsilon }-Y_{r}^{\varepsilon }%
\big\rangle\medskip%
\end{array}%
\end{equation*}%
for $s=t_{i+1}\wedge \tau $,$r=t_{i}\wedge \tau $, where $%
t=t_{0}<t_{1}<t_{2}<...<t\wedge \tau $ and $t_{i+1}-t_{i}=\dfrac{1}{n}$,
then summing up over $i$, and passing to the limit as $n\rightarrow \infty $%
, we deduce%
\begin{equation*}
\begin{array}{r}
\displaystyle e^{\lambda \tau +\mu A_{\tau }}\varphi _{\varepsilon }\left(
\xi \right) \geq e^{\lambda \left( t\wedge \tau \right) +\mu A_{t\wedge \tau
}}\varphi _{\varepsilon }\left( Y_{t\wedge \tau }^{\varepsilon }\right)
+\int_{t\wedge \tau }^{\tau }e^{\lambda s+\mu A_{s}}\big\langle\nabla
\varphi _{\varepsilon }\left( Y_{s}^{\varepsilon }\right)
,dY_{s}^{\varepsilon }\big\rangle \\
\displaystyle+\int_{t\wedge \tau }^{\tau }\varphi _{\varepsilon }\left(
Y_{s}^{\varepsilon }\right) d(e^{\lambda s+\mu A_{s}})%
\end{array}%
\end{equation*}%
We have the similar inequalities for function $\psi _{\varepsilon }$

If we summing and we use equation (\ref{ecapr}), we infer that for all $%
t\geq 0:$\newline
$%
\begin{array}{c}
\displaystyle e^{\lambda \left( t\wedge \tau \right) +\mu A_{t\wedge \tau }}%
\big(\varphi _{\varepsilon }\left( Y_{t\wedge \tau }^{\varepsilon }\right)
+\psi _{\varepsilon }\left( Y_{t\wedge \tau }^{\varepsilon }\right) \big)%
+\int_{t\wedge \tau }^{\tau }e^{\lambda s+\mu A_{s}}\big|\nabla \psi
_{\varepsilon }\left( Y_{s}^{\varepsilon }\right) \big|^{2}dA_{s}\medskip%
\end{array}%
$\newline
$%
\begin{array}{c}
\displaystyle\quad +\int_{t\wedge \tau }^{\tau }e^{\lambda s+\mu A_{s}}\big|%
\nabla \varphi _{\varepsilon }\left( Y_{s}^{\varepsilon }\right) \big|%
^{2}ds+\int_{t\wedge \tau }^{\tau }e^{\lambda s+\mu A_{s}}\big(\varphi
_{\varepsilon }\left( Y_{s}^{\varepsilon }\right) +\psi _{\varepsilon
}\left( Y_{s}^{\varepsilon }\right) \big)(\lambda ds+\mu dA_{s})\medskip%
\end{array}%
$\newline
$%
\begin{array}{c}
\displaystyle\quad +\int_{t\wedge \tau }^{\tau }e^{\lambda s+\mu A_{s}}%
\big\langle\nabla \varphi _{\varepsilon }\left( Y_{s}^{\varepsilon }\right)
,\nabla \psi _{\varepsilon }\left( Y_{s}^{\varepsilon }\right) \big\rangle%
\left( ds+dA_{s}\right) \medskip%
\end{array}%
$\newline
$%
\begin{array}{c}
\displaystyle\leq e^{\lambda \tau +\mu A_{\tau }}\big(\varphi _{\varepsilon
}\left( \xi \right) +\psi _{\varepsilon }\left( \xi \right) \big)%
+\int_{t\wedge \tau }^{\tau }e^{\lambda s+\mu A_{s}}\big\langle\nabla
\varphi _{\varepsilon }\left( Y_{s}^{\varepsilon }\right) ,F\left(
s,Y_{s}^{\varepsilon },Z_{s}^{\varepsilon }\right) \big\rangle ds\medskip%
\end{array}%
$\newline
$%
\begin{array}{c}
\displaystyle\quad +\int_{t\wedge \tau }^{\tau }e^{\lambda s+\mu A_{s}}%
\big\langle\nabla \varphi _{\varepsilon }\left( Y_{s}^{\varepsilon }\right)
,G\left( s,Y_{s}^{\varepsilon }\right) \big\rangle dA_{s}\vspace*{2mm}%
\end{array}%
$\newline
$%
\begin{array}{c}
\displaystyle\quad +\int_{t\wedge \tau }^{\tau }e^{\lambda s+\mu A_{s}}%
\big\langle\nabla \psi _{\varepsilon }\left( Y_{s}^{\varepsilon }\right)
,F\left( s,Y_{s}^{\varepsilon },Z_{s}^{\varepsilon }\right) \big\rangle %
ds\medskip%
\end{array}%
$\newline
$%
\begin{array}{c}
\displaystyle\quad +\int_{t\wedge \tau }^{\tau }e^{\lambda s+\mu A_{s}}%
\big\langle\nabla \psi _{\varepsilon }\left( Y_{s}^{\varepsilon }\right)
,G\left( s,Y_{s}^{\varepsilon }\right) \big\rangle dA_{s}\medskip%
\end{array}%
$\newline
$%
\begin{array}{c}
\displaystyle\quad -\int_{t\wedge \tau }^{\tau }e^{\lambda s+\mu A_{s}}%
\big\langle\nabla \varphi _{\varepsilon }\left( Y_{s}^{\varepsilon }\right)
+\nabla \psi _{\varepsilon }\left( Y_{s}^{\varepsilon }\right)
,Z_{s}^{\varepsilon }dW_{s}\big\rangle\medskip%
\end{array}%
$\newline
The result follows by combining this with (\ref{ineg1}), assumptions (\ref%
{comp}) and the following inequalities\newline
$%
\begin{array}{c}
\dfrac{1}{2\varepsilon }\big|y-J_{\varepsilon }\left( y\right) \big|^{2}\leq
\varphi _{\varepsilon }\left( y\right) \;,\;\dfrac{1}{2\varepsilon }\big|y-%
\hat{J}_{\varepsilon }\left( y\right) \big|^{2}\leq \psi _{\varepsilon
}\left( y\right) \medskip%
\end{array}%
$\newline
$%
\begin{array}{c}
\varphi \big(J_{\varepsilon }\left( y\right) \big)\leq \varphi _{\varepsilon
}\left( y\right) \;,\;\psi \big(\hat{J}_{\varepsilon }\left( y\right) \big)%
\leq \psi _{\varepsilon }\left( y\right) \medskip%
\end{array}%
$\newline
$%
\begin{array}{c}
\varphi _{\varepsilon }\left( \xi \right) \leq \varphi \left( \xi \right)
\;,\;\psi _{\varepsilon }\left( \xi \right) \leq \psi \left( \xi \right)
\medskip%
\end{array}%
$\newline
$%
\begin{array}{c}
\big\langle\nabla \varphi _{\varepsilon }\left( y\right) ,F\left(
s,y,z\right) \big\rangle\leq \dfrac{1}{4}\big|\nabla \varphi _{\varepsilon
}\left( y\right) \big|^{2}+3\big(\eta _{s}^{2}+L^{2}\left\vert y\right\vert
^{2}+L^{2}\left\vert \left\vert z\right\vert \right\vert ^{2}\big)\medskip%
\end{array}%
$\newline
$%
\begin{array}{c}
\big\langle\nabla \psi _{\varepsilon }\left( y\right) ,G\left( s,y\right) %
\big\rangle\leq \dfrac{1}{4}\left\vert \nabla \psi _{\varepsilon }\left(
y\right) \right\vert ^{2}+2\big(\gamma _{s}^{2}+L^{2}\left\vert y\right\vert
^{2}\big)\vspace*{2mm}%
\end{array}%
$\newline
$%
\begin{array}{l}
\big\langle\nabla \psi _{\varepsilon }\left( y\right) ,F\left( s,y,z\right) %
\big\rangle\leq \big\langle\nabla \varphi _{\varepsilon }\left( y\right)
,F\left( s,y,z\right) \big\rangle^{+}\vspace*{2mm} \\
\;\;\;\;\;\;\;\;\;\;\;\;\;\;\;\;\;\;\;\;\;\;\;\;\;\;\;\;\;\;\;\;\;\;\leq
\dfrac{1}{4}\big|\nabla \varphi _{\varepsilon }\left( y\right) \big|^{2}+3%
\big(\eta _{s}^{2}+L^{2}\left\vert y\right\vert ^{2}+L^{2}\left\vert
\left\vert z\right\vert \right\vert ^{2}\big)\medskip%
\end{array}%
$\newline
$%
\begin{array}{l}
\big\langle\nabla \varphi _{\varepsilon }\left( y\right) ,G\left( s,y\right) %
\big\rangle\leq \big\langle\nabla \psi _{\varepsilon }\left( y\right)
,G\left( s,y\right) \big\rangle^{+}\vspace*{2mm} \\
\;\;\;\;\;\;\;\;\;\;\;\;\;\;\;\;\;\;\;\;\;\;\;\;\;\;\;\;\;\;\;\;\;\leq
\dfrac{1}{4}\big|\nabla \psi _{\varepsilon }\left( y\right) \big|^{2}+2\big(%
\gamma _{s}^{2}+L^{2}\left\vert y\right\vert ^{2}\big)\vspace*{2mm}%
\end{array}%
$\hfill
\end{proof}

\begin{proposition}
\label{p4}Let assumptions (\ref{ipF})-(\ref{comp}) be satisfied. Then%
\begin{equation}
\begin{array}{l}
\displaystyle\mathbb{E}\int_{0}^{\tau}\text{ }e^{^{\lambda t+\mu A_{t}}}\big(%
\left\vert Y_{s}^{\varepsilon}-Y_{s}^{\delta}\right\vert ^{2}\left(
ds+dA_{s}\right) +\left\Vert Z_{s}^{\varepsilon}-Z_{s}^{\delta}\right\Vert
^{2}ds\big)\medskip \\
\text{ \ \ \ \ \ \ \ \ \ \ \ \ \ \ \ \ \ \ \ \ \ }\displaystyle+\mathbb{E}%
\underset{0\leq t\leq\tau}{\sup}e^{\lambda t+\mu A_{t}}\left\vert
Y_{t}^{\varepsilon}-Y_{t}^{\delta}\right\vert ^{2}\leq C\left( \varepsilon
+\delta\right) ~M\left( \tau\right)%
\end{array}
\label{ineg3}
\end{equation}
\end{proposition}

\begin{proof}
By It\^{o}'s formula\newline
$%
\begin{array}{c}
\displaystyle e^{\lambda \left( t\wedge \tau \right) +\mu A_{t\wedge \tau
}}\left\vert Y_{t\wedge \tau }^{\varepsilon }-Y_{t\wedge \tau }^{\delta
}\right\vert ^{2}+\int_{t\wedge \tau }^{\tau }e^{\lambda s+\mu
A_{s}}\left\vert Y_{s}^{\varepsilon }-Y_{s}^{\delta }\right\vert ^{2}\left(
\lambda ds+\mu dA_{s}\right) \medskip%
\end{array}%
$\newline
$%
\begin{array}{c}
\displaystyle\quad +2\int_{t\wedge \tau }^{\tau }e^{\lambda s+\mu
A_{s}}\left\langle Y_{s}^{\varepsilon }-Y_{s}^{\delta },\nabla \varphi
_{\varepsilon }\left( Y_{s}^{\varepsilon }\right) -\nabla \varphi _{\delta
}(Y_{s}^{\delta })\right\rangle ds\medskip%
\end{array}%
$\newline
$%
\begin{array}{c}
\displaystyle\quad +2\int_{t\wedge \tau }^{\tau }e^{\lambda s+\mu
A_{s}}\left\langle Y_{s}^{\varepsilon }-Y_{s}^{\delta },\nabla \psi
_{\varepsilon }\left( Y_{s}^{\varepsilon }\right) -\nabla \psi _{\delta
}(Y_{s}^{\delta })\right\rangle dA_{s}=\medskip%
\end{array}%
$\newline
$%
\begin{array}{c}
\displaystyle=2\int_{t\wedge \tau }^{\tau }e^{\lambda s+\mu
A_{s}}\left\langle Y_{s}^{\varepsilon }-Y_{s}^{\delta },F\left(
s,Y_{s}^{\varepsilon },Z_{s}^{\varepsilon }\right) -F(s,Y_{s}^{\delta
},Z_{s}^{\delta })\right\rangle ds\medskip%
\end{array}%
$\newline
$%
\begin{array}{c}
\displaystyle\quad +2\int_{t\wedge \tau }^{\tau }e^{\lambda s+\mu
A_{s}}\left\langle Y_{s}^{\varepsilon }-Y_{s}^{\delta },G\left(
s,Y_{s}^{\varepsilon }\right) -G(s,Y_{s}^{\delta })\right\rangle
dA_{s}\medskip%
\end{array}%
$\newline
$%
\begin{array}{c}
\displaystyle\quad -\int_{t\wedge \tau }^{\tau }e^{\lambda s+\mu
A_{s}}\left\Vert Z_{s}^{\varepsilon }-Z_{s}^{\delta }\right\Vert
^{2}ds-2\int_{t\wedge \tau }^{\tau }e^{\lambda s+\mu A_{s}}\left\langle
Y_{s}^{\varepsilon }-Y_{s}^{\delta },(Z_{s}^{\varepsilon }-Z_{s}^{\delta
})dW_{s}\right\rangle%
\end{array}%
$\newline
We have moreover,
\begin{equation*}
\begin{array}{l}
2\left\langle Y_{s}^{\varepsilon }-Y_{s}^{\delta },F\left(
s,Y_{s}^{\varepsilon },Z_{s}^{\varepsilon }\right) -F(s,Y_{s}^{\delta
},Z_{s}^{\delta })\right\rangle \vspace*{2mm} \\
\;\;\;\;\;\;\;\;\;\;\;\;\;\;\;\;\;\;\;\;\;\;\;\;\;\;\;\;\leq \left( 2\alpha
+2L^{2}\right) \left\vert Y_{s}^{\varepsilon }-Y_{s}^{\delta }\right\vert
^{2}+\dfrac{1}{2}\left\Vert Z_{s}^{\varepsilon }-Z_{s}^{\delta }\right\Vert
^{2}%
\end{array}%
\end{equation*}%
\begin{equation*}
2\left\langle Y_{s}^{\varepsilon }-Y_{s}^{\delta },G\left(
s,Y_{s}^{\varepsilon }\right) -G(s,Y_{s}^{\delta })\right\rangle \leq 2\beta
\left\vert Y_{s}^{\varepsilon }-Y_{s}^{\delta }\right\vert ^{2}
\end{equation*}%
But from definition of $\varphi _{\varepsilon }$ and the monotonicity of
operator $\partial \varphi $ we have%
\begin{equation*}
\begin{array}{l}
0\leq \big\langle\nabla \varphi _{\varepsilon }\left( Y_{s}^{\varepsilon
}\right) -\nabla \varphi _{\delta }(Y_{s}^{\delta }),J_{\varepsilon }\left(
Y_{s}^{\varepsilon }\right) -J_{\delta }(Y_{s}^{\delta })\big\rangle\vspace*{%
2mm} \\
\;\;\;=\big\langle\nabla \varphi _{\varepsilon }\left( Y_{s}^{\varepsilon
}\right) -\nabla \varphi _{\delta }(Y_{s}^{\delta }),Y_{s}^{\varepsilon
}-Y_{s}^{\delta }\big\rangle-\varepsilon \big|\nabla \varphi _{\varepsilon
}\left( Y_{s}^{\varepsilon }\right) \big|^{2}-\delta \big|\nabla \varphi
_{\delta }(Y_{s}^{\delta })\big|^{2}\medskip \\
\;\;\;\;\;\ \;+\left( \varepsilon +\delta \right) \big\langle\nabla \varphi
_{\varepsilon }\left( Y_{s}^{\varepsilon }\right) ,\nabla \varphi _{\delta
}(Y_{s}^{\delta })\big\rangle%
\end{array}%
\end{equation*}%
then%
\begin{equation*}
\big\langle\nabla \varphi _{\varepsilon }\left( Y_{s}^{\varepsilon }\right)
-\nabla \varphi _{\delta }(Y_{s}^{\delta }),Y_{s}^{\varepsilon
}-Y_{s}^{\delta }\big\rangle\geq -\left( \varepsilon +\delta \right) %
\big\langle\nabla \varphi _{\varepsilon }\left( Y_{s}^{\varepsilon }\right)
,\nabla \varphi _{\delta }(Y_{s}^{\delta })\big\rangle
\end{equation*}%
and in the same manner%
\begin{equation*}
\big\langle\nabla \psi _{\varepsilon }\left( Y_{s}^{\varepsilon }\right)
-\nabla \psi _{\delta }(Y_{s}^{\delta }),Y_{s}^{\varepsilon }-Y_{s}^{\delta }%
\big\rangle\geq -\left( \varepsilon +\delta \right) \big\langle\nabla \psi
_{\varepsilon }\left( Y_{s}^{\varepsilon }\right) ,\nabla \psi _{\delta
}(Y_{s}^{\delta })\big\rangle
\end{equation*}%
and consequently%
\begin{equation}
\begin{array}{l}
\displaystyle e^{\lambda \left( t\wedge \tau \right) +\mu A_{t\wedge \tau
}}|Y_{t\wedge \tau }^{\varepsilon }-Y_{t\wedge \tau }^{\delta
}|^{2}+\int_{t\wedge \tau }^{\tau }e^{\lambda s+\mu
A_{s}}|Y_{s}^{\varepsilon }-Y_{s}^{\delta }|^{2}(\lambda -2\alpha
-2L^{2})ds\medskip \\
\displaystyle+\int_{t\wedge \tau }^{\tau }e^{\lambda s+\mu
A_{s}}|Y_{s}^{\varepsilon }-Y_{s}^{\delta }|^{2}\left( \mu -2\beta \right)
dA_{s}+\dfrac{1}{2}\int_{t\wedge \tau }^{\tau }e^{\lambda s+\mu
A_{s}}||Z_{s}^{\varepsilon }-Z_{s}^{\delta }||^{2}ds\medskip \\
\displaystyle\leq 2\left( \varepsilon +\delta \right) \int_{t\wedge \tau
}^{\tau }e^{\lambda s+\mu A_{s}}\Big [\big\langle\nabla \varphi
_{\varepsilon }\left( Y_{s}^{\varepsilon }\right) ,\nabla \varphi _{\delta
}(Y_{s}^{\delta })\big\rangle ds \\
\displaystyle+\big\langle\nabla \psi _{\varepsilon }\left(
Y_{s}^{\varepsilon }\right) ,\nabla \psi _{\delta }(Y_{s}^{\delta })%
\big\rangle dA_{s}\Big ]-2\int_{t\wedge \tau }^{\tau }e^{\lambda s+\mu A_{s}}%
\big\langle Y_{s}^{\varepsilon }-Y_{s}^{\delta },(Z_{s}^{\varepsilon
}-Z_{s}^{\delta })dW_{s}\big\rangle%
\end{array}
\label{ineg4}
\end{equation}%
Now, from (\ref{ineg2}-a)%
\begin{equation*}
\begin{array}{l}
\displaystyle2\left( \varepsilon +\delta \right) \mathbb{E}\int_{t\wedge
\tau }^{\tau }e^{\lambda s+\mu A_{s}}\Big [\big\langle\nabla \varphi
_{\varepsilon }\left( Y_{s}^{\varepsilon }\right) ,\nabla \varphi _{\delta
}(Y_{s}^{\delta })\big\rangle ds\medskip \\
\;\ \ \ \ \ \ \ \ \ \ \ \ \ \ \ \ \ \ \ \ \ \ \ \ \ \ +\big\langle\nabla
\psi _{\varepsilon }\left( Y_{s}^{\varepsilon }\right) ,\nabla \psi _{\delta
}(Y_{s}^{\delta })\big\rangle dA_{s}\Big ]\leq C~\left( \varepsilon +\delta
\right) ~M\left( \tau \right)%
\end{array}%
\end{equation*}%
and clearly by standard calculus inequality (\ref{ineg3}) follows.\hfill $%
\medskip $
\end{proof}

We give now the proof of Theorem \ref{t1}

\begin{proof}
Uniqueness is a consequence of Proposition \ref{p1}. The existence of
solution $\left( Y,Z,U,V\right) $ is obtained as limit of $\left(
Y_{s}^{\varepsilon },Z_{s}^{\varepsilon },\nabla \varphi _{\varepsilon
}\left( Y_{s}^{\varepsilon }\right) ,\nabla \psi _{\varepsilon }\left(
Y_{s}^{\varepsilon }\right) \right) $

From Proposition \ref{p4} we have
\begin{equation*}
\left\vert \ \
\begin{array}{l}
\exists\ Y\in\mathcal{S}_{k}^{\lambda,\mu}\cap\mathcal{H}_{k}^{\lambda,\mu
}\cap\mathcal{\tilde{H}}_{k}^{\lambda,\mu},\ \ \exists\ Z\in\mathcal{H}%
_{k\times d}^{\lambda,\mu}\vspace*{2mm} \\
\underset{{\varepsilon}\searrow0}{\lim}Y^{\varepsilon}=Y\text{ \ in \ }%
\mathcal{S}_{k}^{\lambda,\mu}\cap\mathcal{H}_{k}^{\lambda,\mu}\cap\mathcal{%
\tilde{H}}_{k}^{\lambda,\mu},\vspace*{2mm} \\
\underset{{\varepsilon}\searrow0}{\lim}Z^{\varepsilon}=Z\text{ \ in \ }%
\mathcal{H}_{k\times d}^{\lambda,\mu}.%
\end{array}
\right.
\end{equation*}
Also, from (\ref{ineg2}-a) and (\ref{ineg2}-c) we have%
\begin{equation*}
\underset{{\varepsilon}\searrow0}{\lim}J_{\varepsilon}(Y^{\varepsilon })=Y%
\text{ \ in \ }\mathcal{H}_{k}^{\lambda,\mu},\;\;\underset{{\varepsilon }%
\searrow0}{\lim}\hat{J}_{\varepsilon}(Y^{\varepsilon})=Y\text{ \ in \ }%
\mathcal{\tilde{H}}_{k}^{\lambda,\mu}\vspace*{2mm}
\end{equation*}%
\begin{equation*}
\underset{{\varepsilon}\searrow0}{\lim}\mathbb{E}e^{\lambda\theta+\mu
A_{\theta}}\big |J_{\varepsilon}(Y_{\theta}^{\varepsilon})-Y_{\theta }\big |%
^{2}=0,\;\;\underset{{\varepsilon}\searrow0}{\lim}\mathbb{E}%
e^{\lambda\theta+\mu A_{\theta}}\big |\hat{J}_{\varepsilon}(Y_{\theta
}^{\varepsilon})-Y_{\theta}\big |^{2}=0,
\end{equation*}
for any stopping time $\theta,~\ 0\leq\theta\leq\tau$

Using Fatou's Lemma, from (\ref{ineg2}-b), (\ref{ineg2}-d) and the fact that
$\varphi$ is l.s.c. we obtained (\ref{def}-c).

Denoting $U^{\varepsilon }=\nabla \varphi _{\varepsilon }(Y^{\varepsilon })$%
, $V^{\varepsilon }=\nabla \psi _{\varepsilon }(Y^{\varepsilon })$ , from (%
\ref{ineg2}-a) it follows:%
\begin{equation*}
\mathbb{E}\left[ \int_{\theta }^{\tau }e^{\lambda s+\mu A_{s}}\big (%
|U^{\varepsilon }|^{2}ds+|V^{\varepsilon }|^{2}dA_{s}\big )\right] \leq
C~M\left( \tau \right)
\end{equation*}%
Hence there exists $U\in \mathcal{H}_{k}^{\lambda ,\mu }$ and $V\in \mathcal{%
\tilde{H}}_{k}^{\lambda ,\mu }$ such that for a subsequence $\varepsilon
_{n}\searrow 0$%
\begin{equation*}
\begin{array}{c}
U^{\varepsilon _{n}}\rightharpoonup U,\;\text{weakly in Hilbert space }%
\mathcal{H}_{k}^{\lambda ,\mu } \\
V^{\varepsilon _{n}}\rightharpoonup V,\;\text{weakly in Hilbert space }%
\mathcal{\tilde{H}}_{k}^{\lambda ,\mu }%
\end{array}%
\end{equation*}%
and then%
\begin{equation*}
\begin{array}{l}
\displaystyle\mathbb{E}\left[ \int_{\theta }^{\tau }e^{\lambda s+\mu
A_{s}}\left( \left\vert U\right\vert ^{2}ds+\left\vert V\right\vert
^{2}dA_{s}\right) \right] \medskip \\
\displaystyle\leq \liminf_{n\rightarrow \infty }\mathbb{E}\left[
\int_{\theta }^{\tau }e^{\lambda s+\mu A_{s}}\left( |U^{\varepsilon
_{n}}|^{2}ds+|V^{\varepsilon _{n}}|^{2}dA_{s}\right) \right] \leq
C~M_{2}\left( \theta ,\tau \right) .%
\end{array}%
\end{equation*}%
Passing now to $\lim $ in (\ref{ecapr}) we obtain (\ref{def}-e).

Let $u\in\mathcal{H}_{k}^{\lambda,\mu}$, $v\in\mathcal{\tilde{H}}%
_{k}^{\lambda,\mu}.$ Since $\nabla\varphi_{\varepsilon}(Y_{t}^{\varepsilon})%
\in{\partial\varphi}\big (J_{\varepsilon}(Y_{t}^{\varepsilon})\big )$ and $%
\nabla\psi_{\varepsilon}(Y_{t}^{\varepsilon})\in{\partial\psi}\big (\hat {J}%
_{\varepsilon}(Y_{t}^{\varepsilon})\big ),$ $\forall t\geq0,$ then as signed
measures on $\Omega\times\left[ 0,\tau\right] $%
\begin{equation*}
\begin{array}{r}
{e^{\lambda s+\mu A_{s}}}\big\langle U_{s}^{\varepsilon},u_{s}-J_{%
\varepsilon }\left( {Y_{s}^{\varepsilon}}\right) \big\rangle~\mathbb{P}%
\left( d\omega\right) \otimes ds+{e^{\lambda s+\mu A_{s}}}\varphi \big (%
J_{\varepsilon}\left( {Y_{s}^{\varepsilon}}\right) \big )~\mathbb{P}\left(
d\omega\right) \otimes ds\medskip \\
\leq{e^{\lambda s+\mu A_{s}}}\varphi(u_{s})~\mathbb{P}\left( d\omega\right)
\otimes ds%
\end{array}%
\end{equation*}
and%
\begin{equation*}
\begin{array}{l}
{e^{\lambda s+\mu A_{s}}}\big\langle V_{s}^{\varepsilon},v_{s}-\hat {J}%
_{\varepsilon}\left( {Y_{s}^{\varepsilon}}\right) \big\rangle~\mathbb{P}%
\left( d\omega\right) \otimes A\left( \omega,ds\right) \vspace*{2mm} \\
+{e^{\lambda s+\mu A_{s}}}\psi\big (\hat{J}_{\varepsilon}\left( {%
Y_{s}^{\varepsilon}}\right) \big )~\mathbb{P}\left( d\omega\right) \otimes
A\left( \omega,ds\right) \medskip \\
\leq{e^{\lambda s+\mu A_{s}}}\psi(v_{s})~\mathbb{P}\left( d\omega\right)
\otimes A\left( \omega,ds\right) .%
\end{array}%
\end{equation*}
Passing to $\lim\inf$ in these last two inequalities we obtain (\ref{def}%
-d).\ The proof is complete.\hfill
\end{proof}

\section{PVI - proof of the existence theorem}

\label{sect3}It follows from result in \cite{li-sz/84} that for each $\left(
t,x\right) \in \mathbb{R}_{+}\times \overline{\mathcal{D}}$ there exists a
unique pair of continuous $\mathcal{F}_{s}^{t}$-p.m.s.p. $%
(X_{s}^{t,x},A_{s}^{t,x})_{s\geq 0}$, with values in $\overline{\mathcal{D}}%
\times \mathbb{R}_{+}$, solution of the reflected stochastic differential
equation
\begin{equation}
\left\{
\begin{array}{l}
\displaystyle X_{s}^{t,x}=x+\int_{t}^{s\vee
t}b(r,X_{r}^{t,x})dr+\int_{t}^{s\vee t}\sigma
(r,X_{r}^{t,x})dW_{r}-\int_{t}^{s\vee t}\nabla \ell
(X_{r}^{t,x})dA_{r}^{t,x},\;\medskip \\
s\longmapsto A_{s}^{t,x}\text{\ \ is increasing}\medskip \\
\displaystyle A_{s}^{t,x}=\int_{t}^{s\vee t}\mathbf{1}_{\{X_{r}^{t,x}\in
Bd\left( \mathcal{D}\right) \}}dA_{r}^{t,x},%
\end{array}%
\right.  \label{defx}
\end{equation}%
where%
\begin{equation*}
\mathcal{F}_{s}^{t}=\sigma \left\{ W_{r}-W_{t}:t\leq r\leq s\right\} \vee
\mathcal{N}.
\end{equation*}

Since $\overline{\mathcal{D}}$ is a bounded set, then%
\begin{equation}
\underset{s\geq0}{\sup}|X_{s}^{t,x}|\leq M  \label{mgx}
\end{equation}
and with similar calculus as in \cite{pa-zh/98} we have for all $\mu,T,p>0$
there exists a positive constant $C$ such that$\;\forall~t,t^{^{\prime}}\in%
\left[ 0,T\right] ,~x,x^{\prime}\in\overline{\mathcal{D}}:$

\begin{equation}
\mathbb{E}\underset{s\in\left[ 0,T\right] }{\sup}|X_{s}^{t,x}-X_{s}^{t^{%
\prime},x^{\prime}}|^{p}\leq C\big (|x-x^{\prime}|^{p}+|t-t^{\prime }|^{%
\frac{p}{2}}\big ),  \label{contx}
\end{equation}
and
\begin{equation}
\mathbb{E}[e^{\mu A_{T}^{t,x}}]<\infty.  \label{mga}
\end{equation}

Let $T>0$ be arbitrary fixed. Under assumptions (\ref{h1})-(\ref{h6}), it
follows from Theorem \ref{t1} with $\tau $ replaced by $T$ that for each $%
\left( t,x\right) \in \left[ 0,T\right] \times \overline{\mathcal{D}}$ there
exists a unique solution $(Y^{tx},Z^{tx},U^{tx},V^{tx})\;$of p.m.s.p.%
\begin{equation*}
\begin{array}{l}
Y^{tx}\in \mathcal{S}_{1}^{\lambda ,\mu }\cap \mathcal{H}_{1}^{\lambda ,\mu
}\cap \mathcal{\tilde{H}}_{1}^{\lambda ,\mu },\vspace*{2mm} \\
Z^{tx}\in \mathcal{H}_{d}^{\lambda ,\mu },\;U^{tx}\in \mathcal{H}%
_{1}^{\lambda ,\mu },\;V^{tx}\in \mathcal{\tilde{H}}_{1}^{\lambda ,\mu }%
\end{array}%
\end{equation*}%
with $Y_{s}^{t,x}=Y_{t}^{t,x}$, $Z_{s}^{t,x}=0$, $U_{s}^{t,x}=0$, $%
V_{s}^{t,x}=0$, for all $s\in \left[ 0,t\right] $\newline
solution of BSDE:%
\begin{equation*}
\begin{array}{l}
\displaystyle Y_{s}^{t,x}+\int_{s}^{T}U_{r}^{t,x}dr+%
\int_{s}^{T}V_{r}^{t,x}dA_{r}^{t,x}=h\left( X_{T}^{t,x}\right) \medskip \\
\displaystyle+\int_{s}^{T}1_{\left[ t,T\right] }\left( r\right) ~f\left(
r,X_{r}^{t,x},Y_{r}^{t,x},Z_{r}^{t,x}\right) dr\medskip \\
\displaystyle+\int_{s}^{T}1_{\left[ t,T\right] }\left( r\right) g\left(
r,X_{r}^{t,x},Y_{r}^{t,x}\right)
dA_{r}^{t,x}-\int_{s}^{T}Z_{r}^{t,x}dW_{r},\;\;\text{for all }s\in \left[ 0,T%
\right] \;\text{a.s.}%
\end{array}%
\end{equation*}%
such that $\ \left( Y_{s}^{t,x},U_{s}^{t,x}\right) \in \partial \varphi $, $%
\mathbb{P}\left( d\omega \right) \otimes dt$,$\;\;\left(
Y_{s}^{t,x},V_{s}^{t,x}\right) \in \partial \psi $, $\mathbb{P}\left(
d\omega \right) \otimes A\left( \omega ,dt\right) ,\;\;$a.e. on $\ \Omega
\times \left[ t,T\right] .$

We observe that function $f$, $g\;$depends by $\omega $ only via function $%
X^{t,x}$.\medskip

\begin{proposition}
Under assumptions (\ref{h1})-(\ref{h6}), we have%
\begin{equation}
\mathbb{E}\underset{s\in\left[ 0,T\right] }{\sup}e^{\lambda s+\mu
A_{s}}|Y_{s}^{t,x}|^{2}\leq C\left( T\right)  \label{mgy}
\end{equation}
and%
\begin{equation}
\begin{array}{l}
\displaystyle\mathbb{E}\underset{s\in\left[ 0,T\right] }{\sup}e^{\lambda
s+\mu A_{s}}|Y_{s}^{t,x}-Y_{s}^{t^{\prime},x^{\prime}}|^{2}\leq\mathbb{E}%
\Big[e^{\lambda\tau+\mu A_{\tau}}\big |h(X_{T}^{tx})-h(X_{T}^{t^{\prime
}x^{\prime}})\big |^{2}\medskip \\
\displaystyle+\int_{0}^{T}e^{\lambda r+\mu A_{r}}\big |%
1_{[t,T]}(r)f(r,X_{r}^{tx},Y_{r}^{tx},Z_{r}^{tx})-1_{[t^{%
\prime},T]}(r)f(r,X_{r}^{t^{\prime}x^{\prime}},Y_{r}^{tx},Z_{r}^{tx})\big |%
^{2}dr\medskip \\
\displaystyle+\int_{0}^{T}e^{\lambda r+\mu A_{r}}\big |%
1_{[t,T]}(r)g(r,X_{r}^{tx},Y_{r}^{tx})-1_{[t^{\prime},T]}(r)g(r,X_{r}^{t^{%
\prime }x^{\prime}},Y_{r}^{tx})\big |^{2}dA_{r}^{t,x}\Big]%
\end{array}
\label{conty}
\end{equation}
\end{proposition}

\begin{proof}
Inequality (\ref{mgy}) follows from Theorem \ref{t1} using also (\ref{mgx}),
(\ref{mga}). Inequality (\ref{conty}) follows from (\ref{unic}) in
Proposition \ref{p1}.\hfill \medskip
\end{proof}

We define
\begin{equation}
u(t,x)=Y_{t}^{tx},\ \ \ (t,x)\in\lbrack0,T]\times\overline{\mathcal{D}}
\label{defu}
\end{equation}
which is a determinist quantity since $Y_{t}^{tx}$ is $\mathcal{F}%
_{t}^{t}\equiv\mathcal{N}$--measurable.

From Markov property we have%
\begin{equation}
u(s,X_{s}^{tx})=Y_{s}^{tx}  \label{ecu}
\end{equation}

\begin{corollary}
Under assumptions (\ref{h1})-(\ref{h6}), function $u$ satisfies:%
\begin{equation}
\begin{array}{rl}
\left( a\right) & u{(t,x)}\in Dom\left( \varphi \right) ,\ \ \ \forall {(t,x)%
}\in \lbrack 0,T]\times \overline{\mathcal{D}},\vspace*{2mm} \\
\left( b\right) & u{(t,x)}\in Dom\left( \psi \right) \hfill \hfill ,\ \ \
\forall {(t,x)}\in \lbrack 0,T]\times Bd\left( \mathcal{D}\right) ,\vspace*{%
2mm} \\
\left( c\right) & u\in C\big(\lbrack 0,T]\times \overline{\mathcal{D}}\big).%
\end{array}
\label{propru}
\end{equation}
\end{corollary}

\begin{proof}
Using (\ref{def}-c) we have $\varphi \big(u(t,x)\big)=\mathbb{E}\varphi
(Y_{t}^{tx})<+{\infty }$ and similarly for $\psi .$ Hence (\ref{propru}-a,b)
follows. Let $(t_{n},x_{n})\rightarrow (t,x).$ Then%
\begin{equation*}
\begin{array}{l}
\big |u(t_{n},x_{n})-u(t,x)\big |^{2}=\mathbb{E}%
|Y_{t_{n}}^{t_{n}x_{n}}-Y_{t}^{tx}|^{2}\leq 2\mathbb{E}\sup_{s\in \lbrack
0,T]}|Y_{s}^{t_{n}x_{n}}-Y_{s}^{tx}|^{2} \\
\;\;\;\;\;\;\;\;\;\;\;\;\;\;\;\;\;\;\;\;\;\;\;\;\;\;\;\;\;\;\;\;\;\;\;\;\;\;%
\;\;\;\;\;\;\;\;\;\;\;\;\;\;\;\;\;\;\;\;\;\;\;\;\;\;\;\;\;\;+2\mathbb{E}%
|Y_{t_{n}}^{tx}-Y_{t}^{tx}|^{2}%
\end{array}%
\end{equation*}%
Using (\ref{conty}), (\ref{mgx}), (\ref{contx}) and (\ref{mga}) we obtain $%
u(t_{n},x_{n})\rightarrow u(t,x)$ \ as $(t_{n},x_{n})\rightarrow (t,x)$%
.\hfill
\end{proof}

\bigskip

We present now the proof of Theorem \ref{t2} (existence of the viscosity
solutions).\medskip

\begin{proof}
It suffices to show the existence of the solution of PVI (\ref{pvi}) on an
arbitrary fixed interval $\left[ 0,T\right] $. Setting%
\begin{equation*}
\tilde{u}\left( t,x\right) =u\left( T-t,x\right)
\end{equation*}%
then the existence for problem (\ref{pvi}) it is equivalent with existence
for (\ref{pvi'})%
\begin{equation}
\left\{
\begin{array}{l}
\dfrac{\partial \tilde{u}(t,x)}{\partial t}+\mathcal{\tilde{L}}_{t}\tilde{u}%
\left( t,x\right) +\tilde{f}\big(t,x,\tilde{u}(t,x),(\nabla \tilde{u}\sigma
)(t,x)\big)\in {\partial \varphi }\big(\tilde{u}(t,x)\big), \\
\;\;\;\text{\ \ \ }\;\;\;\text{\ \ \ }\;\;\;\text{\ \ \ }\;\;\;\text{\ \ \ }%
\;\;\;\text{\ \ \ }\;\;\;\text{\ \ \ }\;\;\;\text{\ \ \ }\;\;\;\text{\ \ \ }%
\;\;\;\text{\ \ \ }\;\;\;\text{\ \ \ }t\in \left( 0,T\right) ,\;x\in
\mathcal{D},\medskip \\
-\dfrac{\partial \tilde{u}(t,x)}{\partial n}+\tilde{g}\big(t,x,\tilde{u}(t,x)%
\big)\in {\partial \psi }\big(\tilde{u}(t,x)\big),\;\;t\in \left( 0,T\right)
,\;x\in Bd\left( \mathcal{D}\right) ,\vspace*{2mm} \\
\tilde{u}(T,x)=h(x),\;x\in \overline{\mathcal{D}},%
\end{array}%
\right.  \label{pvi'}
\end{equation}%
where%
\begin{align*}
\tilde{f}\left( t,x,u,z\right) & =f\left( T-t,x,u,z\right) ,\;\;\tilde{g}%
\left( t,x,u\right) =g\left( T-t,x,u\right) \\
\tilde{\sigma}\left( t,x\right) & =\sigma \left( T-t,x\right) ,\;\;\tilde{b}%
\left( t,x\right) =b\left( T-t,x\right)
\end{align*}%
and%
\begin{equation*}
\mathcal{\tilde{L}}_{t}v\left( x\right) =\frac{1}{2}\sum_{i,j=1}^{d}(\tilde{%
\sigma}\tilde{\sigma}^{\ast })_{ij}(t,x)\frac{{\partial }^{2}v\left(
x\right) }{{\partial }x_{i}{\partial }x_{j}}+\sum_{i=1}^{d}\tilde{b}_{i}(t,x)%
\frac{{\partial }v\left( x\right) }{{\partial }x_{i}}.
\end{equation*}%
We denote also%
\begin{align*}
& \tilde{V}\left( t,x,p,q,X\right) \overset{def}{=}-p-\frac{1}{2}\mathrm{Tr}%
\big((\tilde{\sigma}\tilde{\sigma}^{\ast })(t,x)X\big)-\big\langle\tilde{b}%
(t,x),q\big\rangle \\
& \;\;\;\;\;\;\;\;\;\;\;\;\;\;\;\;\;\;\;\;\;\;\;\;\;\;\;\;\;-\tilde{f}\big(%
t,x,\tilde{u}(t,x),q\tilde{\sigma}(t,x)\big).
\end{align*}%
In the sequel, for simplicity we keep notations $b,\sigma ,u,f,g,\mathcal{L}%
,V\mathcal{\ }$instead of $\tilde{b},\tilde{\sigma},\tilde{u},\tilde{f},%
\tilde{g},\mathcal{\tilde{L}},\tilde{V}$ and we shall prove that function $u$
defined by (\ref{defu}) is a viscosity solution of parabolic variational
inequality (\ref{pvi'}). We show only that $u$ is a viscosity subsolution of
(\ref{pvi'}) (the supersolution case is similar).

Let $\left( t,x\right) \in \left[ 0,T\right] \times \overline{\mathcal{D}}$
and $(p,q,X)\in \mathcal{P}^{2,+}u(t,x)$.

1. The proof for the case $x\in\mathcal{D}$ is similar of that from \cite%
{pa-ra/98}.

2. Let $x\in Bd\left( \mathcal{D}\right) $. Suppose, contrary to our claim,
that%
\begin{equation*}
\begin{array}{l}
\min\Big\{V\left( t,x,p,q,X\right) +\varphi_{-}^{\prime}\big (u(t,x)\big )%
,\medskip \\
\;\;\;\;\;\;\;\;\;\;\;\;\;\;\;\;\;\text{\ \ \ \ }\;\;\;\big\langle\nabla
\ell\left( x\right) ,q\big\rangle-g\big (t,x,u(t,x)\big )+\psi_{-}^{\prime }%
\big (u(t,x)\big )\Big\}>0%
\end{array}%
\end{equation*}
and we will find a contradiction.

It follows by continuity of $f$, $g$, $u$, $b$, $\sigma $, $\ell $, left
continuity and monotonicity of $\varphi _{-}^{\prime }$ and $\psi
_{-}^{\prime }$ that there exists $\varepsilon >0$, $\delta >0$ such that
for all $\left\vert s-t\right\vert \leq \delta $,\ $\left\vert
y-x\right\vert \leq \delta $,%
\begin{equation}
\begin{array}{l}
-\left( p+\varepsilon \right) -\dfrac{1}{2}\mathrm{Tr}\,\big ((\sigma \sigma
^{\ast })(s,y)\left( X+\varepsilon I\right) \big )-\big\langle %
b(s,y),q+\left( X+\varepsilon I\right) \left( y-x\right) \big\rangle\medskip
\\
-f\big (s,y,u(s,y),\big (q+\left( X+\varepsilon I\right) \left( y-x\right) %
\big )\sigma (s,y)\Big )+\varphi _{-}^{\prime }\big ((u(s,y)\big )>0,\medskip
\\
\;\;\;\;\;\;\;\;\;\;\;\;\;\;\;\;\;\;\;\;\;\;\;\;\;\;\;\;\;\;\;\;\;\;\;\;\;\;%
\;\;\;\;\;\;\;\;\;\;\;\;\;\;\;\;\;\;\;\;\;\;\;\;\;\;\;\;\;\;\;\;\;\;\;\;\;\;%
\;\;\;\;\;\;\;\;\;\;\;\;\;\;\;\;\;\;\;\;\;\;\text{if }x\in D%
\end{array}
\label{pert1}
\end{equation}%
and%
\begin{equation}
\begin{array}{r}
\big\langle\nabla \ell \left( y\right) ,q+\left( X+\varepsilon I\right)
\left( y-x\right) \big\rangle-g\big (s,y,u(s,y)\big )+\psi _{-}^{\prime }%
\big (u(s,y)\big )>0,\medskip \\
\text{if }x\in Bd\left( \mathcal{D}\right)%
\end{array}
\label{pert2}
\end{equation}%
Now since $(p,q,X)\in \mathcal{P}^{2,+}u(t,x)$ there exists $0<\delta
^{\prime }\leq \delta $ such that%
\begin{equation*}
u(s,y)<\hat{u}(s,y),
\end{equation*}%
for all $s\in \left[ 0,T\right] $, $s\neq t$,$\;y\in \overline{\mathcal{D}}$%
, $y\neq x$ \ such that $\left\vert s-t\right\vert \leq \delta ^{\prime }$,\
$\left\vert y-x\right\vert \leq \delta ^{\prime },$ where%
\begin{equation*}
\hat{u}(s,y)=u(t,x)+\left( p+\varepsilon \right) (s-t)+\left\langle
q,y-x\right\rangle +\frac{1}{2}\big\langle\left( X+\varepsilon I\right)
(y-x),y-x\big\rangle
\end{equation*}%
Let
\begin{equation*}
\nu \overset{def}{=}\inf \left\{ s>t:\;|X_{s}^{t,x}-x|\geq \delta ^{\prime
}\right\}
\end{equation*}%
We note that
\begin{equation*}
(\bar{Y}_{s}^{t,x},\bar{Z}_{s}^{t,x})=\left( Y_{s}^{t,x},Z_{s}^{t,x}\right)
,\;t\leq s\leq \left( t+\delta ^{\prime }\right) \wedge \nu
\end{equation*}%
solves the BSDE%
\begin{equation*}
\left\{
\begin{array}{l}
\displaystyle\bar{Y}_{s}^{t,x}=u\left( \nu ,X_{\nu }^{t,x}\right)
+\int_{s}^{\nu }\big (f(r,X_{r}^{t,x},\bar{Y}_{r}^{t,x},\bar{Z}%
_{r}^{t,x})-U_{r}^{t,x}\big )dr-\int_{s}^{\nu }\bar{Z}_{r}^{t,x}dW_{r}%
\medskip \\
\displaystyle\;\ \ \ \ \ \ \ \ \ \ \ \ \ \ \ \ \ \ \ \ \ \ \ \ \ \ \ \ \ \ \
\ \ \ \ \ \ \ \ \ \ \ \ \ \ \ \ +\int_{s}^{\nu }\big (g(r,X_{r}^{t,x},\bar{Y}%
_{r}^{t,x})-V_{r}^{t,x}\big )dA_{r}^{t,x}\bigskip \\
\left( Y_{s}^{t,x},U_{s}^{t,x}\right) \in \partial \varphi ,\mathbb{P}\left(
d\omega \right) \otimes dt,\;\;\left( Y_{s}^{t,x},V_{s}^{t,x}\right) \in
\partial \psi ,\mathbb{P}\left( d\omega \right) \otimes A\left( \omega
,dt\right) , \\
\;\;\;\text{\ \ \ }\;\;\;\text{\ \ \ }\;\;\;\text{\ \ \ }\;\;\;\text{\ \ \ }%
\;\;\;\text{\ \ \ }\;\;\;\text{\ \ \ }\;\;\;\text{\ \ \ }\;\;\;\text{\ \ \ }%
\;\;\;\text{\ \ \ }\;\;\;\text{\ \ \ }\;\;\;\text{\ \ \ }\;\;\;a.e.\ \text{on%
}\ \ \Omega \times \left[ t,T\right] .%
\end{array}%
\right.
\end{equation*}%
Moreover, it follows from It\^{o}'s formula that
\begin{equation*}
(\hat{Y}_{s}^{t,x},\hat{Z}_{s}^{t,x})=\big (\hat{u}(s,X_{s}^{t,x}),\left(
\nabla \hat{u}\sigma \right) (s,X_{s}^{t,x})\big ),\;t\leq s\leq t+\delta
^{\prime }
\end{equation*}%
satisfies%
\begin{equation*}
\begin{array}{l}
\displaystyle\hat{Y}_{s}^{t,x}=\hat{u}(\nu ,X_{\nu }^{t,x})-\int_{s}^{\nu }%
\Big[\dfrac{\partial \hat{u}(r,X_{r}^{t,x})}{\partial t}+\mathcal{L}_{r}\hat{%
u}(r,X_{r}^{t,x})\Big]dr-\int_{s}^{\nu }\hat{Z}_{r}^{t,x}dW_{r}\medskip \\
\displaystyle\;\ \ \ \ \ \ \ \ \ \ \ \ \ \ \ \ \ \ \ \ \ \ \ \ \ \ \ \ \ \ \
\ \ \ \ \ \ \ \ \ \ \ \ \ \ \ \ \ \ +\int_{s}^{\nu }\big\langle\nabla _{x}%
\hat{u}\left( r,X_{r}^{t,x}\right) ,\nabla \ell (X_{r}^{t,x})\big\rangle %
dA_{r}^{t,x}%
\end{array}%
\end{equation*}%
Let $(\tilde{Y}_{s}^{t,x},\tilde{Z}_{s}^{t,x})=(\hat{Y}_{s}^{t,x}-\bar{Y}%
_{s}^{t,x},\hat{Z}_{s}^{t,x}-\bar{Z}_{s}^{t,x}).$

We have%
\begin{equation*}
\begin{array}{l}
\displaystyle\tilde{Y}_{s}^{t,x}=\big[\hat{u}(\nu,X_{\nu}^{t,x})-u\left(
\nu,X_{\nu}^{t,x}\right) \big]+\int_{s}^{\nu}\Big[-\dfrac{\partial\hat {u}%
(r,X_{r}^{t,x})}{\partial t}-\mathcal{L}_{r}\hat{u}(r,X_{r}^{t,x})\medskip
\\
\displaystyle-f(r,X_{r}^{t,x},\bar{Y}_{r}^{t,x},\bar{Z}%
_{r}^{t,x})+U_{r}^{t,x}\Big]dr-\int_{s}^{\nu}\tilde{Z}_{r}^{t,x}dW_{r}%
\vspace*{2mm} \\
\displaystyle+\int_{s}^{\nu}\Big[\big\langle\nabla_{x}\hat{u}%
(r,X_{r}^{t,x}),\nabla\ell(X_{r}^{t,x})\big\rangle-g(r,X_{r}^{t,x},\bar{Y}%
_{r}^{t,x})+V_{r}^{t,x}\Big]dA_{r}^{t,x}.%
\end{array}%
\end{equation*}
Let
\begin{equation*}
\begin{array}{c}
\overline{\beta}_{s}=\mathcal{L}_{s}\hat{u}(s,X_{s}^{t,x})+f(s,X_{s}^{t,x},%
\bar{Y}_{s}^{t,x},\bar{Z}_{s}^{t,x})\medskip \\
\hat{\beta}_{s}=\mathcal{L}_{s}\hat{u}(s,X_{s}^{t,x})+f(s,X_{s}^{t,x},\bar {Y%
}_{s}^{t,x},\hat{Z}_{s}^{t,x})%
\end{array}%
\end{equation*}
Since $|\hat{\beta}_{s}-\bar{\beta}_{s}|\leq C~|\hat{Z}_{s}^{t,x}-\bar{Z}%
_{s}^{t,x}|$, there exists a bounded $d-$dimensional p.m.s.p. $\left\{
\zeta_{s};0\leq s\leq\nu\right\} $ such that\medskip$\hat{\beta}_{s}-\bar{%
\beta}_{s}=\langle\zeta_{s},\tilde{Z}_{s}^{t,x}\rangle$\medskip

Now
\begin{equation*}
\begin{array}{l}
\displaystyle\tilde{Y}_{s}^{t,x}=\big[\hat{u}(\nu ,X_{\nu }^{t,x})-u(\nu
,X_{\nu }^{t,x})\big]\medskip \\
\;\;\;\;\;\;\;\;\displaystyle+\int_{s}^{\nu }\Big[-\dfrac{\partial \hat{u}%
(r,X_{r}^{t,x})}{\partial t}+\langle \zeta _{r},\tilde{Z}_{r}^{t,x}\rangle -%
\hat{\beta}_{r}+U_{r}^{t,x}\Big]dr\medskip \\
\;\;\;\;\;\;\;\;\displaystyle+\int_{s}^{\nu }\Big[\big\langle\nabla _{x}\hat{%
u}(r,X_{r}^{t,x}),\nabla \ell \left( X_{r}^{t,x}\right) \big\rangle-g\left(
r,X_{r}^{t,x},\bar{Y}_{r}^{t,x}\right) +V_{r}^{t,x}\Big]dA_{r}^{t,x} \\
\;\;\;\;\;\;\;\;\displaystyle-\int_{s}^{\nu }\tilde{Z}_{r}^{t,x}dW_{r}%
\end{array}%
\end{equation*}%
It is easily to see that, for the process
\begin{equation*}
\Gamma _{s}^{t}=\exp \left[ -\frac{1}{2}\int_{t}^{s}\left\vert \zeta
_{r}\right\vert ^{2}dr+\int_{t}^{s}\left\langle \zeta
_{r},dW_{r}\right\rangle \right] ,
\end{equation*}%
we have, from It\^{o}'s formula,%
\begin{equation*}
\Gamma _{s}^{t}=\Gamma _{t}^{t}+\int_{t}^{s}\Gamma _{r}^{t}~\left\langle
\zeta _{r},dW_{r}\right\rangle
\end{equation*}%
and so%
\begin{align*}
d(\tilde{Y}_{s}^{t,x}~\Gamma _{s}^{t})& =\Gamma _{s}^{t}\Big[\dfrac{\partial
\hat{u}(s,X_{s}^{t,x})}{\partial t}+\hat{\beta}_{s}-U_{s}^{t,x}\Big]%
ds+\Gamma _{s}^{t}\langle \tilde{Z}_{r}^{t,x}+\tilde{Y}_{s}^{t,x}\zeta
_{s},dW_{s}\rangle \medskip \\
& +\Gamma _{s}^{t}\Big[\big\langle\nabla _{x}\hat{u}\left(
r,X_{r}^{t,x}\right) ,\nabla \ell (X_{r}^{t,x})\big\rangle+g(s,X_{s}^{t,x},%
\bar{Y}_{s}^{t,x})-V_{s}^{t,x}\Big]dA_{s}^{t,x}
\end{align*}%
Then%
\begin{equation}
\begin{array}{l}
\tilde{Y}_{t}^{t,x}=\mathbb{E}\Big[\Gamma _{\nu }^{t}~\big(\hat{u}(\nu
,X_{\nu }^{t,x})-u(\nu ,X_{\nu }^{t,x})\big)\Big]\vspace*{2mm} \\
\;\;\;\;\displaystyle-\mathbb{E}\bigg[\int_{t}^{\nu }\Gamma _{r}^{t}\Big[%
\dfrac{\partial \hat{u}(r,X_{r}^{t,x})}{\partial t}+\hat{\beta}%
_{r}-U_{r}^{t,x}\Big]dr\bigg]\medskip \\
\;\;\;\;\displaystyle-\mathbb{E}\bigg[\int_{t}^{\nu }\Gamma _{r}^{t}\Big (-%
\big\langle\nabla _{x}\hat{u}\left( r,X_{r}^{t,x}\right) ,\nabla \ell \left(
X_{r}^{t,x}\right) \big\rangle+g(r,X_{r}^{t,x},\bar{Y}_{r}^{t,x})-V_{r}^{t,x}%
\Big )dA_{r}^{t,x}\bigg]%
\end{array}
\label{contr}
\end{equation}%
We first note that $\left( Y_{t},U_{t}\right) \in \partial \varphi $,$%
\;\;\left( Y_{t},V_{t}\right) \in \partial \psi $ implies that%
\begin{equation*}
\varphi _{-}^{\prime }\big (u(s,X_{s}^{t,x})\big )ds\leq U_{s}^{t,x}ds,\;\ \
\psi _{-}^{\prime }\big (u(s,X_{s}^{t,x})\big )dA_{s}^{t,x}\leq
V_{s}^{t,x}dA_{s}^{t,x}
\end{equation*}%
Moreover, the choice of $\delta ^{\prime }$ and $\nu $ implies that
\begin{equation*}
u(\nu ,X_{\nu }^{t,x})<\hat{u}(\nu ,X_{\nu }^{t,x})
\end{equation*}%
From (\ref{pert1}) and (\ref{pert2}) it follows that
\begin{equation*}
-\left( p+\varepsilon \right) -\hat{\beta}_{s}+\varphi _{-}^{\prime }\big (%
u(s,X_{s}^{t,x})\big )>0,\;\text{if }x\in \mathcal{D}
\end{equation*}%
and%
\begin{equation*}
\frac{\partial \hat{u}\left( s,X_{s}^{t,x}\right) }{\partial n}%
-g(s,X_{s}^{t,x},\bar{Y}_{s}^{t,x})+\psi _{-}^{\prime }\big (u(s,X_{s}^{t,x})%
\big )>0,\;\text{if }x\in Bd\left( \mathcal{D}\right)
\end{equation*}%
All these inequalities and equation (\ref{contr}) imply that $\tilde{Y}%
_{t}^{t,x}>0$ and equivalent%
\begin{equation*}
\hat{u}\left( t,x\right) >u\left( t,x\right) ,
\end{equation*}%
which is a contradiction with the definition of $\hat{u}.$ Hence we have%
\begin{equation*}
\min \Big\{V\left( t,x,p,q,X\right) +\varphi _{-}^{\prime }\left(
u(t,x)\right) ,\big\langle\nabla \ell \left( x\right) ,q\big\rangle-g\big (%
t,x,u(t,x)\big )+\psi _{-}^{\prime }\big (u(t,x)\big )\Big\}\leq 0
\end{equation*}%
This proves that $u$ is a viscosity subsolution of (\ref{pvi'}). Symmetric
arguments show that $u$ is also a supersolution; hence $u$ is a viscosity
solution of PVI (\ref{pvi'}).\hfill
\end{proof}

\begin{remark}
If $b,\sigma ,f$ and $g$ do not depend on $t$ then we have a directly a
representation formula for the viscosity solution $u$ of PVI (\ref{pvi}):%
\begin{equation*}
u\left( t,x\right) =Y_{0}^{0,x;t}
\end{equation*}%
where $(Y_{s}^{0,x;t},Z_{s}^{0,x;t},U_{s}^{0,x;t},V_{s}^{0,x;t})_{0\leq
s\leq t}$ \ solution of BSVI%
\begin{equation*}
\begin{array}{l}
\displaystyle Y_{s}^{0,x;t}+\int_{s}^{t}U_{r}^{0,x;t}dr+%
\int_{s}^{t}V_{r}^{0,x;t}dA_{r}^{0,x}=h(X_{t}^{0,x})+%
\int_{s}^{t}f(X_{r}^{t,x},Y_{r}^{0,x;t},Z_{r}^{0,x;t})dr\medskip \\
\;\;\;\;\;\;\;\;\;\;\;\;\;\;\;\;\;\;\;\;\;\;\;\;\;\;\;\displaystyle%
+\int_{s}^{t}g(X_{r}^{0,x},Y_{r}^{0,x;t})dA_{r}^{0,x}-%
\int_{s}^{t}Z_{r}^{0,x;t}dW_{r},\;\;\text{for all }s\in \left[ 0,T\right] \;%
\text{a.s.}%
\end{array}%
\end{equation*}%
and $(X_{s}^{0,x},A_{s}^{0,x})_{0\leq s\leq t}$ solves SDE%
\begin{equation*}
\left\{
\begin{array}{l}
\displaystyle X_{s}^{0,x}=x+\int_{0}^{s}b(X_{r}^{0,x})dr+\int_{0}^{s}\sigma
(X_{r}^{0,x})dW_{r}-\int_{0}^{s}\nabla \ell
(X_{r}^{0,x})dA_{r}^{0,x},\;\medskip \\
s\longmapsto A_{s}^{0,x}\text{\ \ is increasing}\medskip \\
\displaystyle A_{s}^{0,x}=\int_{0}^{s}\mathbf{1}_{\{X_{r}^{0,x}\in Bd\left(
\mathcal{D}\right) \}}dA_{r}^{0,x}\;.%
\end{array}%
\right.
\end{equation*}%
$\medskip $
\end{remark}

\begin{corollary}
We have%
\begin{equation*}
u{(t,x)}\in Dom\left( \partial\varphi\right) ,\ \ \ \forall{(t,x)}\in
\lbrack0,T]\times\mathcal{D}
\end{equation*}
\end{corollary}

\begin{proof}
Let ${(t,x)}$ be fixed. We have two cases:

1) \ $Dom\left( {\partial}\varphi\right) =Dom\left( \varphi\right) ,$ and
so, from (\ref{propru}-a), ${u(t,x)}\in Dom\left( {\partial}\varphi\right) $.

2) $Dom\left( {\partial}\varphi\right) \neq Dom\left( \varphi\right) $. Let $%
b\in Dom\,\varphi\setminus Dom\,(\partial\varphi)$.

Then $b=\sup (Dom\,\varphi )$ or $b=\inf Dom\,\varphi .$ If $b=\sup
(Dom\,\varphi )$ and $u(t,x)=b$, then $(0,0,0)\in \mathcal{P}^{2,+}u(t,x)$
since
\begin{equation*}
u(s,y)\leq {u(t,x)}+o\big (|s-t|+|y-x|^{2}\big )
\end{equation*}%
and from (\ref{subs}) it follows $\varphi _{-}^{\prime }(b)=\varphi
_{-}^{\prime }\big (u(t,x)\big )<\infty $ and consequently $b\in
Dom\,(\partial \varphi )$; a contradiction which shows that $u(t,x)<b.$
Similarly for $b=\inf (Dom\,\varphi ).$\hfill
\end{proof}

\end{document}